\documentclass[12pt,oneside,english,reqno]{amsart}
\usepackage[T1]{fontenc}
\usepackage[latin9]{inputenc}
\usepackage[letterpaper]{geometry}
\usepackage{babel}

\usepackage{textcomp}
\usepackage{amsthm}
\usepackage{amstext}
\usepackage{setspace}
\usepackage{amssymb}
\usepackage[all]{xy}
\usepackage[unicode=true,
 bookmarks=false,
 breaklinks=false,pdfborder={0 0 1},backref=false,colorlinks=false]
 {hyperref}
\hypersetup{pdftitle={Covariant representations of subproduct systems},
 pdfauthor={Ami Viselter}}

\makeatletter
\theoremstyle{plain}
\newtheorem{thm}{Theorem}[section]
  \theoremstyle{definition}
  \newtheorem{defn}[thm]{Definition}
  \theoremstyle{plain}
  \newtheorem{lem}[thm]{Lemma}
  \theoremstyle{remark}
  \newtheorem{rem}[thm]{Remark}
 \theoremstyle{definition}
  \newtheorem{example}[thm]{Example}
  \theoremstyle{plain}
  \newtheorem{prop}[thm]{Proposition}
  \theoremstyle{plain}
  \newtheorem{cor}[thm]{Corollary}
  \theoremstyle{remark}
  \newtheorem*{rem*}{Remark}
  \theoremstyle{plain}
  \newtheorem*{thm*}{Theorem}
  \theoremstyle{plain}
  \newtheorem*{cor*}{Corollary}

\usepackage{mathrsfs}

\numberwithin{equation}{section}
\DeclareMathOperator{\Img}{Im}
\DeclareMathOperator{\linspan}{span}
\DeclareMathOperator*{\slim}{s-lim}

\makeatother

\begin{document}
\global\long\def\e{\varepsilon}
\global\long\def\N{\mathbb{N}}
\global\long\def\Z{\mathbb{Z}}
\global\long\def\Q{\mathbb{Q}}
\global\long\def\R{\mathbb{R}}
\global\long\def\C{\mathbb{C}}
\global\long\def\norm#1{\left\Vert #1\right\Vert }
\global\long\def\newmacroname{}
\global\long\def\bra#1{\left(#1\right)}
\global\long\def\Bra#1{\left[#1\right]}

\global\long\def\H{\mathcal{H}}
\global\long\def\a{\alpha}
\global\long\def\be{\beta}
\global\long\def\l{\lambda}

\global\long\def\Sg{\mathfrak{S}}

\global\long\def\tensor{\otimes}

\global\long\def\A{\forall}
\global\long\def\Aalg{\mathscr A}
\global\long\def\B{\mathscr B}

\global\long\def\D{\mathcal{D}}
\global\long\def\F{\mathcal{F}}
\global\long\def\L{\mathcal{L}}
\global\long\def\M{\mathscr M}
\global\long\def\T{\mathcal{T}}
\global\long\def\K{\mathcal{K}}
\global\long\def\U{\mathcal{U}}
\global\long\def\I{\mathcal{I}}
\global\long\def\J{\mathcal{J}}

\title{Covariant representations of subproduct systems}

\author{Ami Viselter}

\address{Department of Mathematics, Technion - Israel Institute of Technology,
32000 Haifa, Israel}

\email{viselter@tx.technion.ac.il}
\begin{abstract}
A celebrated theorem of Pimsner states that a covariant representation
$T$ of a $C^{*}$-correspondence $E$ extends to a $C^{*}$-representation
of the Toeplitz algebra of $E$ if and only if $T$ is isometric.
This paper is mainly concerned with finding conditions for a covariant
representation of a \emph{subproduct system} to extend to a $C^{*}$-representation
of the Toeplitz algebra. This framework is much more general than
the former. We are able to find sufficient conditions, and show that
in important special cases, they are also necessary. Further results
include the universality of the tensor algebra, dilations of completely
contractive covariant representations, Wold decompositions and von
Neumann inequalities.
\end{abstract}
\maketitle
\tableofcontents{}

\section*{Introduction}

This paper treats the interplay between operator algebras and covariant
representations associated with subproduct systems, introduced and
studied recently in \cite{Subproduct_systems_2009}. This topic is
the intersection of two, not completely disjoint, paths of research.
The first is the investigation of operator algebras corresponding
to product systems. The second concerns with special classes of row
contractions, such as commuting or $q$-commuting or, more generally,
{}``constrained'' row contractions. Both paths, briefly described
hereinafter, greatly generalize much of the Sz.-Nagy--Foia\c{s} harmonic
analysis of Hilbert space contractions (\cite{Nagy_Foias_book}) by
considering families of operators, namely covariant representations,
satisfying some characteristic conditions. In each of these contexts,
certain operator algebras ($C^{*}$ and non-selfadjoint) were proved
to be universal with respect to the features of the covariant representations.

Row contractions have been the subject of extensive study for several
decades (e.g. \cite{Frazho_1982,Bunce_1984,Frazho_1984,Popescu_isom_dialations_1989,Popescu_unit_ball_1991,Popescu_func_calc_1995,Popescu_noncomm_disk_alg_1996,Popescu_noncomm_joint_dil_1998,Popescu_univ_oper_alg_1998,Popescu_Poisson_1999}),
yielding results on isometric dilations, the Wold decomposition, representations
of related operator algebras, the von Neumann inequality and non-commutative
functional calculus, to name a few, with an affluence of applications.
Motivated by the celebrated paper of Pimsner \cite{Pimsner}, Muhly
and Solel generalized in \cite{Tensor_algebras} a large portion of
this theory to the framework of covariant representations indexed
by Hilbert modules, stimulating a variety of deep follow-ups (see
\cite{wold_decomposition,Quotient_Tensor_Algebras,Quantum_Markov_Processes,Hardy algebras 1,Poisson kernel}).

We sketch out some fundamental results from \cite{Pimsner,Tensor_algebras},
which are the basis for this paper, referring to §\ref{sec:Preliminaries-and-notations}
for more details. Suppose that $E$ is a $C^{*}$-correspondence,
that is, a Hilbert $C^{*}$-module over a $C^{*}$-algebra $\M$ equipped
with a left $\M$-module structure implemented by a {*}-homomorphism
from $\M$ to $\L(E)$. The (full) Fock space is $\F(E):=\M\oplus E\oplus E^{\tensor2}\oplus\cdots$.
Given $\zeta\in E$, denote by $S_{\zeta}$ the creation (shift) operator
in $\L(\F(E))$ defined by $E^{\tensor n}\ni\eta\mapsto\zeta\tensor\eta$.
The $C^{*}$-algebra generated by the operators $\left\{ S_{\zeta}:\zeta\in E\right\} $
is called the \emph{Toeplitz algebra}, and is denoted by $\T(E)$,
and the non-selfadjoint operator algebra generated by these operators
is called the \emph{tensor algebra}, and is denoted by $\T_{+}(E)$.
The importance of the Toeplitz algebra lies in the fact that it is
universal in the following sense. Fix a Hilbert space $\H$. A pair
$(T,\sigma)$ is called a \emph{covariant representation} of $E$
on $\H$ if $\sigma$ is a $C^{*}$-representation of $\M$ on $\H$
and $T(\cdot)$ is a bimodule map from $E$ to $B(\H)$. If \begin{equation}
(\A\zeta,\eta\in E)\quad\quad T(\zeta)^{*}T(\eta)=\sigma(\left\langle \zeta,\eta\right\rangle ),\label{eq:intro_isometric_covariant_rep}\end{equation}
then the covariant representation $(T,\sigma)$ is said to be \emph{isometric},
and if $T(\cdot)$ is completely contractive, then $(T,\sigma)$ is
said to be \emph{completely contractive}. (For example, when $\M=\C$
and $E=\C^{d}$, there is a bijection between completely contractive,
covariant representations of $E$ on $\H$ and row contractions of
length $d$ over $\H$ given by $T\mapsto(T(e_{1}),\ldots,T(e_{d}))$;
now \eqref{eq:intro_isometric_covariant_rep} is equivalent to $T(e_{1}),\ldots,T(e_{d})$
being isometries with orthogonal ranges.) It is proved in \cite[Theorem 3.4]{Pimsner}
(with minor differences) that there exists a $C^{*}$-representation
of $\T(E)$ on $\H$, mapping $S_{\zeta}$ to $T(\zeta)$, if and
only if $(T,\sigma)$ is isometric. In other words, the Toeplitz algebra
is the \emph{universal} $C^{*}$-algebra generated by an isometric
covariant representation of $E$.

Concurrently, a theory was sought to match row contractions subject
to restrictions. For example, row contractions consisting of commuting
(\cite{Ando,Drury,Arveson_subalgebras_3,Popescu_Poisson_1999,Bhat_Bhattacharyya_Dey_2004})
or $q$-commuting (\cite{Arias_Popescu_I,Bhat_Bhattacharyya_q_model,Dey})
operators were examined. The general case of constrained row contractions
followed (\cite{Popescu_noncomm_varieties,Popescu_noncomm_varieties_II}),
and proved to have many applications.

Subproduct systems provide the means to unify these two theories.
Assume that $X=\left(X(n)\right)_{n\in\Z_{+}}$ is a family of Hilbert
modules over $\M=X(0)$ such that $X(n+m)$ is an orthogonally complementable
sub-bimodule of $X(n)\tensor X(m)$ for all $n,m\in\Z_{+}$. This
is a subproduct system in the {}``standard'' form. Setting $E:=X(1)$
we have $X(n)\subseteq E^{\tensor n}$ for all $n$. Letting $p_{n}$
denote the orthogonal projection of $E^{\tensor n}$ onto $X(n)$,
the $X$-shifts are defined over the $X$-Fock space $\F_{X}:=\M\oplus E\oplus X(2)\oplus X(3)\oplus\ldots\subseteq\F(E)$
by $S_{n}^{X}(\zeta)\eta:=p_{n+m}(\zeta\tensor\eta)$ for $\zeta\in X(n),\eta\in X(m)$.
The Toeplitz algebra $\T(X)$ and the tensor algebra $\T_{+}(X)$
are now defined to be the $C^{*}$- and non-selfadjoint algebras,
respectively, generated by the $X$-shifts in $\L(\F_{X})$. Covariant
representations $T=\left(T_{n}\right)_{n\in\Z_{+}}$ of $X$ on $\H$
are defined is a suitable manner where $T_{n}:X(n)\to B(\H)$ for
all $n$. Now, for instance, \cite{Pimsner,Tensor_algebras} correspond
to \emph{product} systems ($X(n)=E^{\tensor n}$, so that actually
$X(n+m)=X(n)\tensor X(m)$ for all $n,m$); and there is a bijection
between row contractions of $d$ \emph{commuting} operators and completely
contractive, covariant representations of the \emph{symmetric} subproduct
system defined by $X(n):=(\C^{d})^{\circledS n}$, $n\in\Z_{+}$.

The main goal of this paper is to generalize the above-mentioned theorem
of Pimsner to this setting. That is, if $T$ is a covariant representation
of a subproduct system $X$ on $\H$, \emph{when is there a $C^{*}$-representation
of $\T(X)$ on $\H$ mapping $S_{n}^{X}(\zeta)$ to $T_{n}(\zeta)$}?
If this holds, $T$ is said to \emph{extend to a $C^{*}$-representation}.
The difficulty starts with the fact that the $X$-shifts and their
adjoints do not satisfy any {}``isometricity'' relation in the spirit
of \eqref{eq:intro_isometric_covariant_rep} (as seen even in the
simple example of the symmetric subproduct system). An essential step
in the proof of \cite[Theorem 3.4]{Pimsner} employs \eqref{eq:intro_isometric_covariant_rep}
to reduce every composition of shifts and their adjoints to an operator
of the form $S_{\zeta_{1}}\cdots S_{\zeta_{n}}S_{\eta_{m}}^{*}\cdots S_{\eta_{1}}^{*}$.
Such computation is evidently not possible in the subproduct systems
case, and new ideas will be utilized to establish our results.

The structure of the paper is the following. After presenting some
preliminaries in §\ref{sec:Preliminaries-and-notations}, we introduce
in §\ref{sec:The-Poisson-kernel} the Poisson kernel suitable for
our context. It is then used to prove the universality of the tensor
algebra $\T_{+}(X)$---every completely contractive, covariant representation
extends to a completely contractive representation of $\T_{+}(X)$---and
to prove a dilation theorem for completely contractive, covariant
representations. We next address the question of $C^{*}$-representability
of covariant representations, which is the main part of the work.
In view of our dilation theorem, we divide the problem into two: \emph{pure}
covariant representations are handled in §\ref{sec:Pure-and-relatively_isometric},
while \emph{fully coisometric} covariant representations are dealt
with in §\ref{sec:Fully-coisometric}. In both cases, sufficient conditions
are found for $C^{*}$-representability. Although we are not able
to prove them necessary generally, we demonstrate that this is not
far from being true by providing examples of important special cases
in which equivalence holds. The last section §\ref{sec:Conclusions-and-examples}
is devoted to deriving some conclusions from the results of the paper
and giving more examples.

\section{Preliminaries and notations\label{sec:Preliminaries-and-notations}}

\subsection{Basics}

Throughout this paper, $\H$ denotes an arbitrary (complex) Hilbert
space, and $B(\H)$ denotes the Banach space of bounded operators
over $\H$.

We assume that the reader is familiar with the basic theory of Hilbert
$C^{*}$-modules, which can be found, e.g., in \cite[Ch.~1-4]{Lance}.
As is customary, the action of the $C^{*}$-algebra on the Hilbert
$C^{*}$-module is on the right. The Banach space of all adjointable
(bounded) module maps between two Hilbert $C^{*}$-modules $E,F$
over the same algebra is denoted by $\L(E,F)$, and $\L(E)$ stands
for the $C^{*}$-algebra $\L(E,E)$. A Hilbert $W^{*}$-module is
a Hilbert $C^{*}$-module over a von Neumann algebra that is self-dual
(see \cite[\S 3]{Paschke_1973}, \cite[\S 2.5]{Hilbert_modules_R}).
In this case, $\L(E)$ is a $W^{*}$-algebra. We use the notation
$(\cdot,\cdot)$ for the inner product in Hilbert spaces and $\left\langle \cdot,\cdot\right\rangle $
for the rigging in Hilbert modules.

A Hilbert $C^{*}$-module $E$ over a $C^{*}$-algebra $\M$ is called
a $C^{*}$\emph{-correspondence} in case it is also equipped with
a left $\M$-module structure, implemented by a {*}-homomorphism $\varphi:\M\to\L(E)$;
that is, $a\cdot\zeta:=\varphi(a)\zeta$ for $a\in\M$, $\zeta\in E$.
In particular, $\left\langle a\cdot\zeta,\eta\right\rangle =\left\langle \zeta,a^{*}\cdot\eta\right\rangle $.
$E$ is called a $W^{*}$\emph{-correspondence} when $E$ is a Hilbert
$W^{*}$-module (in particular, $\M$ is a von Neumann algebra) and
$\varphi$ is normal.

Let $\M,\mathscr N$ be $C^{*}$-algebras, and let $E,F$ be Hilbert
$C^{*}$-modules over $\M,\mathscr N$ respectively. Suppose also
that $\sigma:\M\to\L(F)$ is a {*}-homomorphism. The \emph{(interior)
tensor product} of $E$ and $F$, $E\tensor_{\sigma}F$, is the Hilbert
$C^{*}$-module over $\mathscr N$ that contains the algebraic tensor
product of $E$ and $F$, balanced by $\sigma$ (i.e., $(\zeta a)\tensor\eta=\zeta\tensor(\sigma(a)\eta)$),
as a dense subspace, with the rigging defined by\begin{alignat*}{1}
(\A\zeta_{1},\zeta_{2}\in E,\eta_{1},\eta_{2}\in F)\qquad & \left\langle \zeta_{1}\tensor\eta_{1},\zeta_{2}\tensor\eta_{2}\right\rangle _{E\tensor_{\sigma}F}:=\left\langle \eta_{1},\sigma(\left\langle \zeta_{1},\zeta_{2}\right\rangle )\eta_{2}\right\rangle _{F}.\end{alignat*}
 Two special cases are worth mentioning. First, if $F=\H$ is a Hilbert
space, i.e. $\mathscr N=\C$, then $E\tensor_{\sigma}\H$ is also
a Hilbert space. Second, if $E,F$ are $C^{*}$-correspondences over
the same algebra $\M$ (and $\varphi_{F}$ implements the left action
of $\M$ on $F$), so is $E\tensor F:=E\tensor_{\varphi_{F}}F$.

Assume that $E_{1},E_{2}$ are Hilbert $C^{*}$-modules over $\M$
and $F_{1},F_{2}$ are Hilbert $C^{*}$-modules over $\mathscr N$.
Suppose also that $\sigma_{j}:\M\to\L(F_{j})$, $j=1,2$, are {*}-homomorphisms.
If $T\in\L(E_{1},E_{2}),S\in\L(F_{1},F_{2})$ and $S$ is an $\M$-module
map, that is, $S\sigma_{1}(a)=\sigma_{2}(a)S$, then there exists
a unique operator $T\tensor S\in\L(E_{1}\tensor_{\sigma_{1}}F_{1},E_{2}\tensor_{\sigma_{2}}F_{2})$
satisfying $(T\tensor S)(\zeta\tensor\eta)=(T\zeta)\tensor(S\eta)$
for $\zeta\in E_{1},\eta\in F_{1}$. Moreover, $\norm{T\tensor S}\leq\norm T\cdot\norm S$.
Tensor products of Hilbert $W^{*}$-modules are discussed more thoroughly
in Remark \ref{rem:self_dual_completion}.

For a general (not necessarily selfadjoint) operator algebra $\mathscr A$,
denote by $\mathscr A^{1}$ its trivial unitization $\mathscr A\oplus\C I$
if it is not already unital (otherwise, set $\mathscr A^{1}:=\mathscr A$).
It should be emphasized that the $C^{*}$-algebras we consider are
not necessarily unital. However, all operator algebra representations,
and particularly $C^{*}$-algebra representations, are assumed nondegenerate
by convention.

The majority of the theory developed in this paper is valid in both
the $C^{*}$-algebra setting and the $W^{*}$ (von Neumann)-algebra
setting. To avoid burdensome repetitions in statements, we shall use
the letter $A$ where either $C$ or $W$ may be used, yielding, e.g.,
an $A^{*}$-algebra, a Hilbert $A^{*}$-module, etc. \emph{Unless
specifically told otherwise, all definitions, propositions, etc. are
valid in both contexts}.

\subsection{Covariant representations of Hilbert modules}

Fix an $A^{*}$-correspondence $E$ over an $A^{*}$-algebra $\M$.
\begin{defn}
[{\cite[Definition 2.11]{Tensor_algebras}, \cite[Definition 2.15]{Quantum_Markov_Processes}}]\label{def:c_c_representation_correspondence}A
pair $(T,\sigma)$ is called a \emph{covariant representation} of
$E$ on $\H$ if:
\begin{enumerate}
\item $\sigma$ is a $C^{*}$-representation of $\M$ on $\H$.
\item $T$ is a linear mapping from $E$ to $B(\H)$.
\item $T$ is a bimodule map with respect to $\sigma$, i.e., $T(\zeta a)=T(\zeta)\sigma(a)$
and $T(a\zeta)=\sigma(a)T(\zeta)$ for all $\zeta\in E$ and $a\in\M$.
\end{enumerate}
In the $W^{*}$-setting we require, in addition, that $\sigma$ be
normal and that $T$ be continuous with respect to the $\sigma$-topology
of $E$ and the ultraweak topology on $B(\H)$ (see \cite[pp.~201-202]{Baillet-Denizeau-Havet}).

A covariant representation $(T,\sigma)$ is\emph{ }said to be \emph{completely
contractive} if $T$ is a completely contractive map, when $E$ is
viewed as a subspace of the {}``linking algebra'' of $\M$ and $E$
(see \cite[pp.~398-399]{Tensor_algebras}). It is said to be \emph{isometric}
when $T(\zeta)^{*}T(\eta)=\sigma(\left\langle \zeta,\eta\right\rangle )$
for all $\zeta,\eta\in E$. An isometric, covariant representation
is necessarily completely contractive.
\end{defn}
Given a completely contractive, covariant representation $(T,\sigma)$
of $E$ on $\H$, define an operator $\widetilde{T}:E\tensor_{\sigma}\H\to\H$
by \begin{equation}
\widetilde{T}(\zeta\tensor h):=T(\zeta)h\label{eq:T_tilde_def}\end{equation}
 for $\zeta\in E$, $h\in\H$. It turns out that $\widetilde{T}$
is a well-defined contraction, and that if the left action of $\M$
on $E$ is given by the mapping $\varphi:\M\to\L(E)$, then \begin{equation}
(\A a\in\M)\quad\quad\widetilde{T}(\varphi(a)\tensor I_{\H})=\sigma(a)\widetilde{T}.\label{eq:covariance_condition_tilde}\end{equation}
Moreover, we have the following useful lemma, for which we assume
that $\sigma$ is a \emph{fixed} (normal, in the $W^{*}$-setting)
representation of $\M$ on $\H$.
\begin{lem}
[{\cite[Lemma 3.5]{Tensor_algebras}; \cite[Lemma 2.16]{Quantum_Markov_Processes}}; see also the succeeding remark]\label{lem:bijection__c_c_c_rep__tilde}The
mapping $(T,\sigma)\mapsto\widetilde{T}$ defined above is a bijection
between the set of all completely contractive, covariant representations
$(T,\sigma)$ of $E$ on $\H$ and the set of all contractions $\widetilde{T}:E\tensor_{\sigma}\H\to\H$
satisfying \eqref{eq:covariance_condition_tilde}. Additionally, $(T,\sigma)$
is isometric if and only if $\widetilde{T}$ is an isometric operator.\end{lem}
\begin{rem}
\label{rem:sigma_T_normality}It is a consequence of the lemma that,
in the $W^{*}$-setting, the ultraweak continuity of $T$ is automatic
when $\sigma$ is normal.
\end{rem}

\subsection{General subproduct systems}

In this subsection we let $\Sg$ denote a general Abelian monoid (a
semigroup with identity element $0$).
\begin{defn}
[{\cite[Definition 1.1]{Subproduct_systems_2009}}]\label{def:subproduct_system}Let
$\M$ be an $A^{*}$-algebra. A family $X=\left(X(s)\right)_{s\in\Sg}$
of $A^{*}$-correspondences over $\mathcal{\M}$ is called a \emph{subproduct
system over }$\mathcal{\M}$ if the following conditions hold:
\begin{enumerate}
\item \label{enu:subproduct_system_X_0}$X(0)=\mathcal{\M}$.
\item \label{enu:subproduct_system_U}For all $s,t\in\Sg$ there exists
a coisometric, adjointable bimodule mapping \[
U_{s,t}:X(s)\tensor X(t)\to X(s+t),\]
so that:

\begin{enumerate}
\item \label{enu:product_system_U_0}The maps $U_{s,0}$ and $U_{0,s}$
are given by the right and left actions of $\mathcal{\M}$ on $X(s)$,
respectively; that is, $U_{s,0}(\zeta\tensor a)=\zeta a$ and $U_{0,s}(a\tensor\zeta)=a\zeta$
for all $a\in\mathcal{\M},$ $\zeta\in X(s)$.
\item The following associativity condition holds for all $s,t,r\in\Sg$:\[
U_{s+t,r}\left(U_{s,t}\tensor I_{X(r)}\right)=U_{s,t+r}\left(I_{X(s)}\tensor U_{t,r}\right).\]

\end{enumerate}
\end{enumerate}
\end{defn}
\begin{rem}
\label{rem:essential_left_modules}Every Hilbert $C^{*}$-module $E$
over $\M$ is identified with $E\tensor\M$ (via $\zeta\tensor a\mapsto\zeta a$).
However, when $E$ is a $C^{*}$-correspondence, $\M\tensor E$ may
only be embedded in $E$ (via $a\tensor\zeta\mapsto a\zeta$). Equality
holds if and only if $E$ is \emph{essential} as a left $\M$-module,
that is, the {*}-homomorphism $\varphi$ implementing the left multiplication
is nondegenerate: $\varphi(\M)E$ is dense in $E$. If $\M$ is unital
(e.g. in the $W^{*}$-setting), this is equivalent to having $\varphi(I_{\M})=I_{\L(E)}$.
Definition \ref{def:subproduct_system} implies that $X(t)$ is essential
for all $t\in\Sg$. \end{rem}
\begin{example}
A \emph{product system} is a subproduct system for which all mappings
$U_{s,t}$ are \emph{unitary}.
\end{example}

\begin{example}
[Subproduct system of Hilbert spaces]Suppose that $\M=\C$ and $X(s)$
is a Hilbert space for all $s\in\Sg$. This important special case
has been recently studied in \cite{Bhat_Mukherjee}.
\end{example}
We refer the reader to \cite{Subproduct_systems_2009} for many other
examples of subproduct systems.
\begin{rem}
\label{rem:self_dual_completion}We make a few important comments
on tensor products and direct sums of $W^{*}$-correspondences and
bounded operators on them, which will be used, sometimes tacitly,
throughout. The following facts and related definitions are taken
from \cite[\S 3]{Paschke_1973} unless stated otherwise. Fix a von
Neumann algebra $\M$. For an arbitrary Hilbert $C^{*}$-module $X$
over $\M$, recall that $X'$, the linear space of bounded module
maps from $X$ to $\M$, is a \emph{self-dual} Hilbert $C^{*}$-module
over $\M$, that is, $X'$ is a Hilbert $W^{*}$-module. We refer
to $X'$ as the {}``self-dual completion'' of $X$. Hilbert $W^{*}$-modules
possess a variety of desirable properties. For example, every bounded
operator (module map) over a Hilbert $W^{*}$-module is automatically
adjointable, and the space of all such operators is a $W^{*}$-algebra.

We mention some results that demonstrate the connection between $X$
and $X'$. Firstly, for an arbitrary Hilbert $C^{*}$-module $X$
over $\M$, let the \emph{$s$-topology} of $X$ (see \cite[p.~202]{Baillet-Denizeau-Havet})
be defined by the family of semi-norms $(p_{\omega})_{\omega\in\M_{*}^{+}}$
when $p_{\omega}$ is defined by $X\ni x\mapsto\omega(\left\langle x,x\right\rangle )^{1/2}$.
Then $X$ is dense in $X'$ relative to the $s$-topology of the latter
(see \cite[Lemma 2.3]{Paschke_1976} and \cite[Proposition 1.4]{Baillet-Denizeau-Havet}).
Secondly, assume that $X,Y$ are Hilbert $C^{*}$-modules over $\M$.
Then every bounded module map $T:X\to Y$ admits a \emph{unique} extension
to a bounded, adjointable operator $T:X'\to Y'$.

Suppose that $E,F$ are $W^{*}$-correspondences over $\M$. By the
{}``tensor product'' $E\tensor F$ we do not mean the usual interior
tensor product of \cite[Ch.~4]{Lance} (which is denoted by $E\tensor_{C^{*}}F$
for the moment), but rather its self-dual completion. There are hence
several essential subtleties that one has to pay attention to. For
example, the subspace $E\tensor_{\text{{alg}}}F:=\linspan\left\{ \zeta\tensor\eta:\zeta\in E,\eta\in F\right\} $
is not necessarily norm-dense in $E\tensor F$, but rather $s$-dense
there. This {}``problem'' is circumvented when defining bounded
operators from $E\tensor F$ to another Hilbert $W^{*}$-module $G$
by the unique extension property mentioned above. As to comparison
of operators, if $T_{1},T_{2}:E\tensor F\to G$ are bounded module
maps that agree on $E\tensor_{\text{{alg}}}F$, then $\left(T_{1}\right)_{|E\tensor_{C^{*}}F}=\left(T_{2}\right)_{|E\tensor_{C^{*}}F}$,
so that $T_{1}=T_{2}$ by virtue of the uniqueness of the extension
of the operators $T_{i}$ from $E\tensor_{C^{*}}F$ to $E\tensor F$.

If $X$ is a Hilbert $C^{*}$-module over $\M$, $\H$ is a Hilbert
space and $\sigma$ is a \emph{normal} representation of $\M$ on
$\H$, then we may form the Hilbert \emph{space} $X'\tensor_{\sigma}\H$
in the usual way. Thus, using the definition of the $s$-topology
and the normality of $\sigma$, the subspace $X\tensor_{\text{{alg}}}\H:=\linspan\left\{ \zeta\tensor h:\zeta\in X,h\in\H\right\} $
is readily seen to be weakly dense in $X'\tensor_{\sigma}\H$. Since
a linear subspace is convex, $X\tensor_{\text{{alg}}}\H$ is actually
norm-dense in $X'\tensor_{\sigma}\H$. Consequently, $X'\tensor_{\sigma}\H=X\tensor_{\sigma}\H$.
This observation is relevant, e.g., when $E,F$ are $W^{*}$-correspondences
and $X=E\tensor_{C^{*}}F$.

Finally, we comment that if $\left(E_{i}\right)_{I}$ is a family
of Hilbert $W^{*}$-modules over $\M$, its direct sum $\bigoplus_{I}E_{i}$
(which is a Hilbert $C^{*}$-module over $\M$) is not necessarily
self-dual. Its self-dual completion is called the \emph{ultraweak
direct sum} of $\left(E_{i}\right)_{I}$, and is also denoted by $\bigoplus_{I}E_{i}$
where there is no chance of confusion. If $E_{i}$ is a $W^{*}$-correspondence
over $\M$ for all $i\in I$ then $\bigoplus_{I}E_{i}$ is a $W^{*}$-correspondence
over $\M$ too.\end{rem}
\begin{defn}
[{\cite[Definition 1.5]{Subproduct_systems_2009}}]\label{def:c_c_c_representation}Assume
that $\mathcal{\M}$ is an $A^{*}$-algebra and $X=\left(X(s)\right)_{s\in\Sg}$
is a subproduct system over $\mathcal{\M}$. Suppose that $T=\left(T_{s}\right)_{s\in\Sg}$
is a family of linear maps $T_{s}:X(s)\to B(\H)$, and denote $\sigma:=T_{0}$.
Then $T$ is called a \emph{completely contractive, covariant representation
}of \emph{$X$} on $\H$ if:
\begin{enumerate}
\item \label{enu:c_c_representation_condition_threads}For each $s\in\Sg$,
the pair $(T_{s},\sigma)$ is a completely contractive, covariant
representation of the $A^{*}$-correspondence $X(s)$ on $\H$ in
the sense of Definition \ref{def:c_c_representation_correspondence}.
\item \label{enu:c_c_representation_condition}For each $s,t\in\Sg$, $\zeta\in X(s)$
and $\eta\in X(t)$,\begin{equation}
T_{s+t}(U_{s,t}(\zeta\tensor\eta))=T_{s}(\zeta)T_{t}(\eta).\label{eq:c_c_representation_condition}\end{equation}

\end{enumerate}
\end{defn}
An alternative formulation of Condition \ref{enu:c_c_representation_condition}
is obtained by using the bijection described in Lemma \ref{lem:bijection__c_c_c_rep__tilde}.
For $s\in\Sg$ we define the contractions $\widetilde{T}_{s}:X(s)\tensor_{\sigma}\H\to\H$
by $\widetilde{T}_{s}(\zeta\tensor h):=T_{s}(\zeta)h$ for all $\zeta\in X(s)$,
$h\in\H$. The equality in \eqref{eq:c_c_representation_condition}
is now equivalent to \begin{equation}
\widetilde{T}_{s+t}(U_{s,t}\tensor I_{\H})=\widetilde{T}_{s}(I_{X(s)}\tensor\widetilde{T}_{t})\label{eq:c_c_representation_condition_tilde}\end{equation}
(see Remark \ref{rem:self_dual_completion} in this connection). Observe
that $I_{X(s)}\tensor\widetilde{T}_{t}$ is a well-defined bounded
operator from $X(s)\tensor(X(t)\tensor_{\sigma}\H)$ to $X(s)\tensor_{\sigma}\H$
only because \eqref{eq:covariance_condition_tilde} is satisfied for
$T=T_{t}$. Operator tensor products of this form will be frequently
used throughout.
\begin{lem}
\label{lem:c_c_c_representation_T_0_tilde_isometry_onto}Under the
assumptions of Definition \ref{def:c_c_c_representation}, the operator
$\widetilde{T}_{0}:\M\tensor_{\sigma}\H\to\H$ is unitary.\end{lem}
\begin{proof}
That $\widetilde{T}_{0}$ is an isometry is a matter of straightforward
calculation. It is onto since $\sigma$ is nondegenerate by Definition
\ref{def:c_c_c_representation}, \ref{enu:c_c_representation_condition_threads}.
\end{proof}
Let a subproduct system $X=\left(X(s)\right)_{s\in\Sg}$ be given.
We describe next the most natural example of a completely contractive,
covariant {}``representation'' of $X$. The \emph{$X$-Fock space}
is the $A^{*}$-correspondence over $\mathcal{\M}$ defined to be
the (ultraweak, in the $W^{*}$-setting) direct sum of $A^{*}$-correspondences
\[
\F_{X}:=\bigoplus_{s\in\Sg}X(s).\]
For $s\in\Sg$, let $S_{s}^{X}:X(s)\to\L(\F_{X})$ be the linear mapping
defined by \[
S_{s}^{X}(\zeta)\eta:=U_{s,t}^{X}(\zeta\tensor\eta)\]
for all $t\in\Sg$, $\eta\in X(t)$. The operators $S_{s}^{X}$, $s\neq0$,
are called the \emph{creation operators} of $\F_{X}$, and the family
$S^{X}:=(S_{s}^{X})_{s\in\Sg}$ is called the \emph{$X$-shift}. When
there is no ambiguity we write $S$ instead of $S^{X}$.

The $X$-Fock space $\F_{X}$ being an $A^{*}$-correspondence but
not necessarily a Hilbert space, the $X$-shift may not be a covariant
representation as is; this is overcome by composing the mappings $S_{s}^{X}$
with a faithful representation of $\L(\F_{X})$ on some Hilbert space.
The requirements of Definition \ref{def:c_c_c_representation} are
satisfied, so the $X$-shift $S^{X}$ is indeed a completely contractive,
covariant {}``representation'' of $\F_{X}$.

From \ref{enu:subproduct_system_X_0} and \ref{enu:product_system_U_0}
of Definition \ref{def:subproduct_system} we infer that $S_{0}^{X}=\varphi_{\infty}$,
where $\varphi_{\infty}(a)$ maps $(b,\zeta_{1},\zeta_{2},\ldots)\in\F_{X}$
to $(ab,a\cdot\zeta_{1},a\cdot\zeta_{2},\ldots)\in\F_{X}$ for all
$a\in\M=X(0)$. Assuming the $W^{*}$-setting, we claim that the {}``representation''
map $\varphi_{\infty}$ is normal, that is, continuous with respect
to the ultraweak topologies of $\M$ and $\L(\F_{X})$. Let $(a_{i})_{I}$
be a bounded net in $\M$ that converges ultraweakly to $a\in\M$.
Then $\left(\varphi_{\infty}(a_{i})\right)_{I}$ is bounded in $\L(\F_{X})$,
thus it converges ultraweakly to some $A\in\L(\F_{X})$ if and only
if $\omega(\left\langle \zeta,\varphi_{\infty}(a_{i})\eta\right\rangle )\to\omega(\left\langle \zeta,A\eta\right\rangle )$
for all $\omega\in\M_{*}$, $k\in\N$ and $\zeta,\eta\in\bigoplus_{n=0}^{k}X(n)\subseteq\F_{X}$
(see the proof of Proposition 3.10 and the preceding paragraph in
\cite{Paschke_1973} and Remark \ref{rem:self_dual_completion}).
But $\omega(\left\langle \zeta,\varphi_{\infty}(a_{i})\eta\right\rangle )=\sum_{n=0}^{k}\omega(\left\langle \zeta_{n},a_{i}\cdot\eta_{n}\right\rangle )$,
and since, for $n\in\Z_{+}$, $X(n)$ is a $W^{*}$-correspondence,
the action of left multiplication $\M\to\L(X(n))$ is normal, and
we deduce that $\omega(\left\langle \zeta,\varphi_{\infty}(a_{i})\eta\right\rangle )\to\omega(\left\langle \zeta,\varphi_{\infty}(a)\eta\right\rangle )$,
hence $\varphi_{\infty}$ is normal.
\begin{lem}
\textup{\label{lem:phi_infinity_homeomorphism}In the}\emph{ $W^{*}$-setting}\textup{\emph{,
}}\textup{the image $\varphi_{\infty}(\M)$} \emph{is a von Neumann
algebra, and} \textup{$\varphi_{\infty}$ is a homeomorphism of $\M$
onto $\varphi_{\infty}(\M)$, both endowed with their ultraweak topologies.}\end{lem}
\begin{proof}
We have just seen that $\varphi_{\infty}:\M\to\L(\F_{X})$ is normal.
Since it is visibly also faithful, all we need is to apply Proposition
7.1.15 and Corollary 7.1.16 of \cite{Kadison_Ringrose_2}.
\end{proof}

\section{The Poisson kernel, representations of the tensor algebra and dilations\label{sec:The-Poisson-kernel}}

Henceforth we focus on subproduct systems relative to the monoid $\Z_{+}$.
This section is dedicated to adapting two theorems to the framework
of subproduct systems. The first, Theorem \ref{thm:c_c_c_representation__to__A_morphism},
implies that every completely contractive, covariant representation
extends to a representation of the tensor algebra. It is a direct
generalization of, e.g., \cite[\S 8]{Subproduct_systems_2009}. The
second, Theorem \ref{thm:dilation_1}, is a dilation theorem for completely
contractive, covariant representations based on \cite[Theorem 2.1]{Popescu_noncomm_varieties}.

The following is an adaptation of \cite[Definition 6.2]{Subproduct_systems_2009}
to the $C^{*}$-setting as well.
\begin{defn}
\label{def:standard_subproduct_system}A family $X=\left(X(n)\right)_{n\in\Z_{+}}$
of $A^{*}$-correspondences over an $A^{*}$-algebra $\M$ is a \emph{standard}
subproduct system if $X(0)=\M$ and for all $n,m\in\Z_{+}$, $X(n+m)$
is an orthogonally complementable sub-bimodule of $X(n)\tensor X(m)$.
Setting $E:=X(1)$, we have $X(n)\subseteq E^{\tensor n}$ (recall
that $E^{\tensor0}$ equals $\M$ by definition). Let $p_{n}\in\L(E^{\tensor n})$
stand for the orthogonal projection of $E^{\tensor n}$ onto $X(n)$,
and define $P$ to be $\bigoplus_{n\in\Z_{+}}p_{n}$, the orthogonal
projection of $\F(E)$ onto $\F_{X}$.
\end{defn}
The definition implies that the projections $\left(p_{n}\right)_{n\in\Z_{+}}$
are bimodule maps, and \[
(\A n,m\in\mathbf{\Z_{+}})\quad\quad p_{n+m}=p_{n+m}(I_{E^{\tensor n}}\tensor p_{m})=p_{n+m}(p_{n}\tensor I_{E^{\tensor m}}).\]
Every standard subproduct system is a subproduct system over $\M$---simply
define $U_{n,m}:=\left(p_{n+m}\right)_{|X(n)\tensor X(m)}$. By virtue
of \cite[Lemma 6.1]{Subproduct_systems_2009}, we lose nothing by
considering only standard subproduct systems. Although this lemma
is stated there for the $W^{*}$-setting alone, its proof holds verbatim
for the $C^{*}$-setting as well. In the $W^{*}$-setting, the assumption
that $X(n+m)$ is orthogonally complementable is superfluous (see
\cite[Proposition 2.5.4]{Hilbert_modules_R}).
\begin{rem}
The identifications of $(X(n)\tensor X(m))\tensor X(k)$ with $X(n)\tensor(X(m)\tensor X(k))$
and of $X(n)\tensor\M$ and $\M\tensor X(n)$ with $X(n)$ (cf.~Remark
\ref{rem:essential_left_modules}) are implicitly employed in the
preceding definition. \end{rem}
\begin{example}
A standard \emph{product} system $X_{E}$ satisfies $X_{E}(n)=E^{\tensor n}$
for all $n\in\Z_{+}$. Here $E$ may be an arbitrary $C^{*}$-correspondence
that is essential as a left $\M$-module. In the $C^{*}$-setting,
$T\mapsto T_{1}$ is a bijection between the completely contractive,
covariant representations of $X_{E}$ and the completely contractive,
covariant representations of $E$.
\end{example}
Note that if $X$ is a standard subproduct system and $P$ is as above,
then $S_{n}^{X}(\zeta)=PS_{n}^{X_{E}}(\zeta)_{|\F_{X}}$ for all $n\in\Z_{+}$,
$\zeta\in X(n)$.
\begin{example}
\label{exa:finite_dim_Hilbert_spaces}Subproduct systems of finite
dimensional Hilbert spaces constitute a very important special case.
They were studied thoroughly in \cite{Arias_Popescu_I,Popescu_noncomm_varieties,Popescu_noncomm_varieties_II},
mainly in the context of the non-commutative analytic Toeplitz algebra,
which is the weak closure of our tensor algebra (see Definition \ref{def:tensor_Toeplitz_algebras}),
called there the non-commutative disc algebra. Nevertheless, in many
respects, the results of the present section can be viewed as generalizing
parts of \cite{Arias_Popescu_I} and \cite[\S 2]{Popescu_noncomm_varieties}.
Particular examples were explored in \cite{Subproduct_systems_2009,Tsirelson_graded_dim_2,Tsirelson_subproduct_dim_2}.
\end{example}
Assume that $X=\left(X(n)\right)_{n\in\Z_{+}}$ is a standard subproduct
system, and $T$ is a completely contractive, covariant representation
of $X$. Then \eqref{eq:c_c_representation_condition} takes the form\begin{equation}
T_{n+m}(p_{n+m}(\zeta\tensor\eta))=T_{n}(\zeta)T_{m}(\eta),\label{eq:c_c_representation_condition_standard}\end{equation}
and \eqref{eq:c_c_representation_condition_tilde} is now \begin{equation}
\widetilde{T}_{n+m}(p_{n+m}\tensor I_{\H})_{|X(n)\tensor X(m)\tensor_{\sigma}\H}=\widetilde{T}_{n}(I_{X(n)}\tensor\widetilde{T}_{m}).\label{eq:c_c_representation_condition_tilde_standard}\end{equation}
Taking adjoints in \eqref{eq:c_c_representation_condition_tilde_standard}
we hence obtain \begin{equation}
\widetilde{T}_{n+m}^{*}=(I_{X(n)}\tensor\widetilde{T}_{m}^{*})\widetilde{T}_{n}^{*}\label{eq:c_c_representation_condition_tilde_standard_adjoints}\end{equation}
(formally, in the left side of \eqref{eq:c_c_representation_condition_tilde_standard_adjoints}
we should have composed $\widetilde{T}_{n+m}^{*}$ with the embedding
of $X(n+m)\tensor_{\sigma}\H$ in $X(n)\tensor X(m)\tensor_{\sigma}\H$,
but we omit it for the sake of convenience). In particular, \begin{equation}
(\A n\in\Z_{+})\quad\quad\widetilde{T}_{n+1}^{*}=(I_{E}\tensor\widetilde{T}_{n}^{*})\widetilde{T}_{1}^{*}=(I_{X(n)}\tensor\widetilde{T}_{1}^{*})\widetilde{T}_{n}^{*}.\label{eq:T_n_1_tilde_adjoint_n}\end{equation}
Iterating \eqref{eq:c_c_representation_condition_tilde_standard_adjoints}
yields the formula\begin{equation}
(\A n\in\N)\quad\quad\widetilde{T}_{n}^{*}=\bigl(I_{X(n-1)}\tensor\widetilde{T}_{1}^{*}\bigr)\bigl(I_{X(n-2)}\tensor\widetilde{T}_{1}^{*}\bigr)\cdots\bigl(I_{E}\tensor\widetilde{T}_{1}^{*}\bigr)\widetilde{T}_{1}^{*}.\label{eq:T_n_tilde_adjoint_iterative}\end{equation}

\begin{rem}
\label{rem:T_n_tilde_T_n_tilde_adj}Equation \eqref{eq:T_n_1_tilde_adjoint_n}
implies that $\bigl\{\widetilde{T}_{n}\widetilde{T}_{n}^{*}\bigr\}_{n\in\Z_{+}}$
is a decreasing sequence of positive contractions, hence $\slim_{n\to\infty}\widetilde{T}_{n}\widetilde{T}_{n}^{*}$
exists (where $\slim$ stands for limit in the strong operator topology).\end{rem}
\begin{defn}
\label{def:pure_ccc_representation}Let $X=\left(X(n)\right)_{n\in\Z_{+}}$
be a standard subproduct system. A completely contractive, covariant
representation $T$ of $X$ is called \emph{pure} if $\slim_{n\to\infty}\widetilde{T}_{n}\widetilde{T}_{n}^{*}=0$.\end{defn}
\begin{example}
\mbox{}
\begin{enumerate}
\item If $\D$ is a Hilbert space and $\pi$ is a $C^{*}$-representation
of $\M$ on $\D$, then the covariant representation induced by $\pi$
is the family $S\tensor I_{\D}:=\left(S_{n}(\cdot)\tensor I_{\D}\right)_{n\in\Z_{+}}$,
which consists of the induced representation of the operators $S_{n}(\zeta)$
on $\F_{X}\tensor_{\pi}\D$. $S\tensor I_{\D}$ is a pure, completely
contractive, covariant representation of $X$ on $\F_{X}\tensor_{\pi}\D$
\item If $\bigl\Vert\widetilde{T}_{1}\bigr\Vert<1$, then \eqref{eq:T_n_tilde_adjoint_iterative}
yields that $\bigl\Vert\widetilde{T}_{n}\widetilde{T}_{n}^{*}\bigr\Vert\leq\bigl\Vert\widetilde{T}_{1}\bigr\Vert^{2n}$
for all $n$, hence $T$ is pure.
\end{enumerate}
\end{example}
\begin{defn}
\label{def:Delta_star}Let $X=\left(X(n)\right)_{n\in\Z_{+}}$ be
a standard subproduct system, and let $T$ be a completely contractive,
covariant representation of $X$. Define \[
\Delta_{*}(T):=(I_{\H}-\widetilde{T}_{1}\widetilde{T}_{1}^{*})^{\frac{1}{2}}\in B(\H),\quad\D:=\overline{\Img\Delta_{*}(T)}.\]
Then $\Delta_{*}(T)$ is clearly a positive contraction, and it is
invertible in case $\bigl\Vert\widetilde{T}_{1}\bigr\Vert<1$. When
the context is understood, we drop the $T$ and simply write $\Delta_{*}$.\end{defn}
\begin{prop}
\label{pro:delta_star_in_sigma_M_comm}Let $X=\left(X(n)\right)_{n\in\Z_{+}}$
be a standard subproduct system, and let $T$ be a completely contractive,
covariant representation of $X$. Then $\Delta_{*}(T)\in\sigma(\M)'$.\end{prop}
\begin{proof}
It is enough to prove that $\widetilde{T}_{1}\widetilde{T}_{1}^{*}$
belongs to $\sigma(\M)'$. Since $T_{1}$ is a completely contractive,
covariant representation of $E=X(1)$, this equality is a result of
\cite[Lemma 3.6]{Tensor_algebras}.
\end{proof}
By the last proposition, we may consider the bounded operator $I_{X(n)}\tensor\Delta_{*}(T)\in B(X(n)\tensor_{\sigma}\H)$,
and, in a similar fashion, the operator $I_{\F_{X}}\tensor\Delta_{*}(T)\in B(\F_{X}\tensor_{\sigma}\H)$.

Furthermore, $\D$ reduces $\sigma(a)$ for all $a\in\M$. Upon denoting
by $\sigma'$ the reduced representation, one can form the tensor
product Hilbert space $X(n)\tensor_{\sigma'}\D$ for each $n\in\Z_{+}$,
as well as $\F_{X}\tensor_{\sigma'}\D$. For the sake of simplicity,
we write $\sigma$ instead of $\sigma'$.

The operator-related Poisson kernel has played an important role since
its introduction (see \cite{Vasilescu,Arveson_subalgebras_3,Popescu_Poisson_1999}).
The following definition should come as no surprise in light of the
corresponding ones in \cite{Popescu_noncomm_varieties,Poisson kernel}.
\begin{defn}
\label{def:poisson_kernel}Let $X=\left(X(n)\right)_{n\in\Z_{+}}$
be a standard subproduct system, and let $T$ be a completely contractive,
covariant representation of $X$. The \emph{Poisson kernel} of $T$
is the operator $K(T):\H\to\F_{X}\tensor_{\sigma}\D$ defined by \begin{equation}
K(T)h:=\bigoplus_{n\in\Z_{+}}(I_{X(n)}\tensor\Delta_{*}(T))\widetilde{T}_{n}^{*}h\label{eq:poisson_kernel_def}\end{equation}
for all $h\in\H$. It is established in the following proposition
that $K(T)$ is well-defined as an element of $B(\H,\F_{X}\tensor_{\sigma}\D)$.\end{defn}
\begin{prop}
\label{pro:Poisson_kernel_isometric}Under the assumptions of Definition
\ref{def:poisson_kernel}, $K(T)$ is a contraction. It is an isometry
if and only if $T$ is pure.\end{prop}
\begin{proof}
Given $n\in\Z_{+}$ and $h\in\H$, we compute

\begin{equation}
\begin{split}\bigl\Vert(I_{X(n)}\tensor\Delta_{*}(T))\widetilde{T}_{n}^{*}h\bigr\Vert^{2} & =\bigl(\widetilde{T}_{n}(I_{X(n)}\tensor\Delta_{*}(T)^{2})\widetilde{T}_{n}^{*}h,h\bigr)\\
 & =\bigl(\widetilde{T}_{n}(I_{X(n)}\tensor(I_{\H}-\widetilde{T}_{1}\widetilde{T}_{1}^{*}))\widetilde{T}_{n}^{*}h,h\bigr)\\
 & =\bigl(\widetilde{T}_{n}\widetilde{T}_{n}^{*}h,h\bigr)-\bigl(\widetilde{T}_{n+1}\widetilde{T}_{n+1}^{*}h,h\bigr),\end{split}
\label{eq:Poisson_kernel_isometric_1}\end{equation}
where the last equality is deduced from \eqref{eq:c_c_representation_condition_tilde_standard}
and \eqref{eq:T_n_1_tilde_adjoint_n}. Thus\[
\sum_{n=0}^{\infty}\bigl\Vert(I_{X(n)}\tensor\Delta_{*}(T))\widetilde{T}_{n}^{*}h\bigr\Vert^{2}=(h,h)-\lim_{n\to\infty}\bigl(\widetilde{T}_{n}\widetilde{T}_{n}^{*}h,h\bigr)\leq(h,h)\]
($\widetilde{T}_{0}\widetilde{T}_{0}^{*}=I_{\H}$ by Lemma \ref{lem:c_c_c_representation_T_0_tilde_isometry_onto}),
and we conclude that $K(T)$ is a well-defined contraction.

From \eqref{eq:poisson_kernel_def} and the boundedness of $K(T)$
one infers that if $\bigoplus_{\Z_{+}}y_{n}=y\in\F_{X}\tensor_{\sigma}\D$
(where $y_{n}\in X(n)\tensor_{\sigma}\D$ for all $n$), then \begin{equation}
K(T)^{*}y=\sum_{n=0}^{\infty}\widetilde{T}_{n}(I_{X(n)}\tensor\Delta_{*}(T))y_{n}.\label{eq:poisson_kernel_adj}\end{equation}
For $h\in\H$ we hence have\begin{equation}
\begin{split}K(T)^{*}K(T)h & =\sum_{n=0}^{\infty}\widetilde{T}_{n}(I_{X(n)}\tensor(I_{\H}-\widetilde{T}_{1}\widetilde{T}_{1}^{*}))\widetilde{T}_{n}^{*}h\\
 & =(I_{\H}-\widetilde{T}_{1}\widetilde{T}_{1}^{*})h+\sum_{n=1}^{\infty}(\widetilde{T}_{n}\widetilde{T}_{n}^{*}-\widetilde{T}_{n+1}\widetilde{T}_{n+1}^{*})h\\
 & =(I_{\H}-\slim_{n\to\infty}\widetilde{T}_{n}\widetilde{T}_{n}^{*})h.\end{split}
\label{eq:K_adj_K}\end{equation}
Therefore $K(T)^{*}K(T)=I_{\H}$ if and only if $T$ is pure (Definition
\ref{def:pure_ccc_representation}).\end{proof}
\begin{prop}
\label{pro:Poisson_kernel_S_T}Let $X=\left(X(n)\right)_{n\in\Z_{+}}$
be a standard subproduct system, and let $T$ be a completely contractive,
covariant representation of $X$. The following equality holds for
all $n\in\Z_{+}$ and $\zeta\in X(n)$:\[
K(T)^{*}\left(S_{n}(\zeta)\tensor I_{\D}\right)=T_{n}(\zeta)K(T)^{*}.\]
\end{prop}
\begin{proof}
We may consider only vectors of the form $\eta\tensor h$ when $\eta\in X(m)$
for some $m\in\Z_{+}$ and $h\in\D$. Indeed we have, owing to \eqref{eq:c_c_representation_condition_tilde_standard}
and \eqref{eq:poisson_kernel_adj}, \[
\begin{split}K(T)^{*}\left(S_{n}(\zeta)\tensor I_{\D}\right)(\eta\tensor h) & =K(T)^{*}(p_{n+m}(\zeta\tensor\eta)\tensor h)\\
 & =\widetilde{T}_{n+m}(p_{n+m}(\zeta\tensor\eta)\tensor\Delta_{*}(T)h)\\
 & =\widetilde{T}_{n}(\zeta\tensor\widetilde{T}_{m}(\eta\tensor\Delta_{*}(T)h))=T_{n}(\zeta)K(T)^{*}(\eta\tensor h).\qedhere\end{split}
\]
\end{proof}
\begin{defn}
\label{def:tensor_Toeplitz_algebras}Let $X=\left(X(n)\right)_{n\in\Z_{+}}$
be a standard subproduct system. We define the \emph{tensor algebra}
$\T_{+}(X)$ of $X$ to be the (non-selfadjoint, norm closed) subalgebra
of $\L(\F_{X})$ generated by $\left\{ S_{n}(\zeta):n\in\Z_{+},\zeta\in X(n)\right\} $.
The \emph{Toeplitz algebra} $\T(X)$ of $X$ is the $C^{*}$-subalgebra
of $\L(\F_{X})$ generated by these operators. We also define $\mathcal{E}(X)$
to be the operator system $\overline{\linspan}(\T_{+}(X)^{1}\T_{+}(X)^{1*})$,
the closure being taken in the norm operator topology.\end{defn}
\begin{rem}
\label{rem:tensor_algebra}Since the $X$-shift is a completely contractive,
covariant {}``representation'', we deduce from \eqref{eq:c_c_representation_condition_standard}
that actually $\T_{+}(X)=\overline{\linspan}\left\{ S_{n}(\zeta):n\in\Z_{+},\zeta\in X(n)\right\} $.
\end{rem}
The next result implies that $\T_{+}(X)$ is the universal operator
algebra generated by a completely contractive, covariant representation
of $X$.
\begin{thm}
\label{thm:c_c_c_representation__to__A_morphism}Let $X=\left(X(n)\right)_{n\in\Z_{+}}$
be a standard subproduct system, and suppose that $T$ is a completely
contractive, covariant representation of $X$. Then there exists a
unique unital, completely positive, completely contractive linear
map $\Psi:\mathcal{E}(X)\to B(\H)$ that satisfies\begin{equation}
(\A n\in\Z_{+},\zeta\in X(n))\quad\quad\Psi(S_{n}(\zeta))=T_{n}(\zeta)\text{ and }\Psi(S_{n}(\zeta)^{*})=T_{n}(\zeta)^{*}\label{eq:c_c_c_representation__to__A_morphism__S_to_T}\end{equation}
 and \begin{equation}
(\A a\in\T_{+}(X),b\in\mathcal{E}(X))\quad\quad\Psi(ab)=\Psi(a)\Psi(b).\label{eq:c_c_c_representation__to__A_morphism__morphism}\end{equation}
\end{thm}
\begin{proof}
We begin by handling the {}``strict'' case in which $\bigl\Vert\widetilde{T}_{1}\bigr\Vert<1$.
Define a linear map $\Psi:\L(\F_{X})\to B(\H)$ by \[
\Psi(a):=K(T)^{*}\left(a\tensor I_{\D}\right)K(T)\]
for all $a\in\L(\F_{X})$. $\Psi$ is plainly completely contractive
and completely positive. Since $T$ is pure, $K(T)$ is isometric
by Proposition \ref{pro:Poisson_kernel_isometric}. Therefore $\Psi$
is unital, and the rest of the stated features are a consequence of
Proposition \ref{pro:Poisson_kernel_S_T}. Indeed, \begin{equation}
\begin{split}\Psi(S_{n}(\zeta)S_{m}(\eta)^{*}) & =K(T)^{*}\left(S_{n}(\zeta)S_{m}(\eta)^{*}\tensor I_{\D}\right)K(T)\\
 & =K(T)^{*}\left(S_{n}(\zeta)\tensor I_{\D}\right)\left(S_{m}(\eta)^{*}\tensor I_{\D}\right)K(T)\\
 & =T_{n}(\zeta)K(T)^{*}K(T)T_{m}(\eta)^{*}=T_{n}(\zeta)T_{m}(\eta)^{*}\end{split}
\label{eq:c_c_c_representation__to__A_morphism__S_to_T__Psi}\end{equation}
for all $n,m\in\Z_{+}$, $\zeta\in X(n)$ and $\eta\in X(m)$ (and
the proof of \eqref{eq:c_c_c_representation__to__A_morphism__S_to_T}
is similar). This, the norm continuity of $\Psi(\cdot)$ and Remark
\ref{rem:tensor_algebra} yield that \eqref{eq:c_c_c_representation__to__A_morphism__morphism}
is satisfied as stated.

We proceed to the general case. Let $T$ be a completely contractive,
covariant representation of $X$. Given $0<r<1$, define a completely
contractive, covariant representation $_{r}T$ of $X$ (on the same
Hilbert space $\H$) by $_{r}T_{n}:=r^{n}T_{n}$ for all $n\in\Z_{+}$.
Then $\bigl\Vert{_{r}}\widetilde{T}_{1}\bigr\Vert\leq r$, so we can
associate with $_{r}T$ a function $_{r}\Psi$ as above. By \eqref{eq:c_c_c_representation__to__A_morphism__S_to_T__Psi}
we now have \[
_{r}\Psi(S_{n}(\zeta)S_{m}(\eta)^{*})={_{r}}T_{n}(\zeta){_{r}}T_{m}(\eta)^{*}=r^{n+m}T_{n}(\zeta)T_{m}(\eta)^{*}\xrightarrow[r\to1^{-}]{}T_{n}(\zeta)T_{m}(\eta)^{*}\]
for all $n,m\in\Z_{+},\zeta\in X(n),\eta\in X(m)$, the limit being
taken in the norm operator topology. Similarly, $_{r}\Psi(S_{n}(\zeta))\xrightarrow[r\to1^{-}]{}T_{n}(\zeta)$
and $_{r}\Psi(S_{n}(\zeta)^{*})\xrightarrow[r\to1^{-}]{}T_{n}(\zeta)^{*}$.
Hence, by Remark \ref{rem:tensor_algebra}, we conclude that $\lim_{r\to1^{-}}({_{r}}\Psi(a))$
exists for all $a$ in a dense subspace of $\mathcal{E}(X)$. The
family $\left\{ _{r}\Psi:0<r<1\right\} $, which consists of completely
contractive linear maps, is uniformly bounded by $1$. Hence, the
limit $\Psi(a):=\lim_{r\to1^{-}}(_{r}\Psi(a))$ exists for all $a\in\mathcal{E}(X)$,
forming a linear operator $\Psi:\mathcal{E}(X)\to B(\H)$, which is
obviously unital, completely contractive and completely positive,
and satisfies \eqref{eq:c_c_c_representation__to__A_morphism__S_to_T}
and \eqref{eq:c_c_c_representation__to__A_morphism__morphism} (from
which it follows that $\Psi$ is unique).
\end{proof}
Restricting $\Psi$ to $\T_{+}(X)$ we get a completely contractive
\emph{representation} of the operator algebra $\T_{+}(X)$. The next
corollary asserts that by this we obtain a bijection. It generalizes
the analogous results for product systems, \cite[Theorem 3.10]{Tensor_algebras}
and \cite[Theorem 2.9]{Hardy algebras 1}. In the $W^{*}$-setting,
a representation of \emph{$\varphi_{\infty}(\M)$ }on a Hilbert space
$\H$ is called\emph{ normal} if it is continuous with respect to
the ultraweak topologies of $\L(\F_{X})$ and $B(\H)$.
\begin{cor}
Let $X$ be a standard subproduct system. Then there exists a bijection
between the completely contractive, covariant representations of $X$
and the completely contractive representations of the tensor algebra
$\T_{+}(X)$ whose restriction to $\varphi_{\infty}(\M)$ is normal
in the $W^{*}$-setting. This bijection is implemented as follows:
a completely contractive, covariant representation $T$ on a Hilbert
space $\H$ is mapped to the (unique) completely contractive representation
$\rho$ of $\T_{+}(X)$ on $\H$ determined by the equality \begin{equation}
\rho(S_{n}(\zeta))=T_{n}(\zeta)\label{eq:tensor_algebra_c_c_representations_S_to_T}\end{equation}
 for all $n\in\Z_{+}$ and $\zeta\in X(n)$. Additionally, each such
representation is completely positive, and it extends to a unital,
completely positive, completely contractive linear map over $\mathcal{E}(X)$
satisfying \eqref{eq:c_c_c_representation__to__A_morphism__S_to_T}
and \eqref{eq:c_c_c_representation__to__A_morphism__morphism}.\end{cor}
\begin{proof}
By Theorem \ref{thm:c_c_c_representation__to__A_morphism}, there
exists a mapping $T\mapsto\rho$ of a completely contractive, covariant
representation $T$ of $X$ to a completely contractive representation
$\rho$ of $\T_{+}(X)$ that satisfies \eqref{eq:tensor_algebra_c_c_representations_S_to_T}.
The restriction of $\rho$ to $\varphi_{\infty}(\M)$ satisfies $\rho(\varphi_{\infty}(a))=\rho(S_{0}(a))=T_{0}(a)=\sigma(a)$
for all $a\in\M$, and therefore it is normal in the $W^{*}$-setting
by virtue of the normality of $\sigma$ and Lemma \ref{lem:phi_infinity_homeomorphism}.
The last part in the statement of the corollary is also clear from
the theorem. This map is injective by \eqref{eq:tensor_algebra_c_c_representations_S_to_T}.

On the other hand, given a completely contractive representation $\rho$
of $\T_{+}(X)$ on a Hilbert space $\H$ whose restriction to $\varphi_{\infty}(\M)$
is normal in the $W^{*}$-setting, define \[
T_{n}(\zeta):=\rho(S_{n}(\zeta))\]
 for all $n\in\Z_{+}$ and $\zeta\in X(n)$. Since the $X$-shift
is a completely contractive, covariant {}``representation'' of $X$,
one can verify that the conditions of Definition \ref{def:c_c_c_representation}
are fulfilled, thus obtaining a completely contractive, covariant
representation $T$ of $X$ on $\H$ (see, in particular, Lemma \ref{lem:phi_infinity_homeomorphism}
and Remark \ref{rem:sigma_T_normality}, and notice that $\sigma:=T_{0}$
is nondegenerate by Remark \ref{rem:essential_left_modules} because
$\rho$ is nondegenerate). Consequently, our mapping $T\mapsto\rho$
is also a surjection, and the proof is complete.
\end{proof}
Our next conclusion is an adaptation of von Neumann's inequality to
our setting. For this we require the notion of a polynomial.
\begin{defn}
Let $X=\left(X(n)\right)_{n\in\Z_{+}}$ be a standard subproduct system.
A \emph{polynomial} over $X$ is a tuple $(\a,\zeta_{0},\zeta_{1},\ldots,\zeta_{M})$
with $\a\in\C$, $M\in\Z_{+}$ and $\zeta_{n}\in X(n)$ for all $0\leq n\leq M$.
If $T$ is a completely contractive, covariant representation of $X$
on $\H$, define $p(T)\in B(\H)$ to be $\a I_{\H}+\sum_{n=0}^{M}T_{n}(\zeta_{n})$.
$p(S)\in\L(\F_{X})$ is defined similarly.\end{defn}
\begin{cor}
Let $X$ be a standard subproduct system, and suppose that $T$ is
a completely contractive, covariant representation of $X$ on $\H$.
If $p_{1},\ldots,p_{t}$ and $q_{1},\ldots,q_{t}$ are polynomials
over $X$, then \[
\biggl\Vert\sum_{i=1}^{t}p_{i}(T)q_{i}(T)^{*}\biggr\Vert_{B(\H)}\leq\biggl\Vert\sum_{i=1}^{t}p_{i}(S)q_{i}(S)^{*}\biggr\Vert_{\L(\F_{X})}.\]
\end{cor}
\begin{proof}
The corollary follows from Theorem \ref{thm:c_c_c_representation__to__A_morphism}
due to $\Psi$ being contractive.
\end{proof}
We move on to our dilation theorem, beginning with several preliminaries.
\begin{defn}
Let $X=\left(X(n)\right)_{n\in\Z_{+}}$ be a standard subproduct system,
and let $T,V$ be completely contractive, covariant representations
of $X$ on the Hilbert spaces $\H,\K$ respectively, where $\H$ is
a subspace of $\K$. We say that $V$ is a \emph{dilation} of $T$
if for all $n\in\Z_{+}$ and $\zeta\in X(n)$, $V_{n}(\zeta)$ leaves
$\H^{\perp}=\K\ominus\H$ invariant and $P_{\H}V_{n}(\zeta)_{|\H}=T_{n}(\zeta)$
($P_{\H}$ denoting the orthogonal projection of $\K$ on $\H$);
equivalently: $V_{n}(\zeta)^{*}$ leaves $\H$ invariant and $\bigl(V_{n}(\zeta)^{*}\bigr)_{|\H}=T_{n}(\zeta)^{*}$.
\end{defn}
This definition is consistent with \cite[Definition 5.4]{Subproduct_systems_2009}
(when $X=Y$) and with the standard one for product systems (e.g.
\cite[Definition 3.1]{Tensor_algebras}, \cite[Definition 2.18]{Quantum_Markov_Processes}).
As is customary, if $\H,\K$ are Hilbert spaces and $W:\H\to\K$ is
isometric, we regard $\H$ as a subspace of $\K$.
\begin{defn}
Let $X=\left(X(n)\right)_{n\in\Z_{+}}$ be a standard subproduct system.
A completely contractive, covariant representation $T$ of $X$ on
$\H$ is said to be\emph{ isometric} or \emph{fully coisometric} if
$\widetilde{T}_{n}$ is isometric or coisometric, respectively, for
all $n\in\Z_{+}$.
\end{defn}
In both cases, it is enough to check for $n=1$ (see \eqref{eq:c_c_representation_condition_tilde_standard},
\eqref{eq:c_c_representation_condition_tilde_standard_adjoints} and
Lemma \ref{lem:c_c_c_representation_T_0_tilde_isometry_onto}). If
$X$ is a product system, then by Lemma \ref{lem:bijection__c_c_c_rep__tilde},
$T$ is isometric if and only if $(T_{1},\sigma)$ is isometric as
a covariant representation of $E$.

As an introduction to Theorem \ref{thm:dilation_1} and the next sections,
we summarize several results on \emph{product} systems.
\begin{thm}
[{\cite[Theorem 3.4]{Pimsner}, \cite[Theorems 2.12, 3.3]{Tensor_algebras}, \cite[Theorem 2.9]{wold_decomposition}, \cite[Theorem 2.18]{Quantum_Markov_Processes}}]\label{thm:isom_prod_sys}Let
$X=\left(X(n)\right)_{n\in\Z_{+}}$ be a standard product system and
let $V$ be a covariant representation of $X$ on $\H$. Then:
\begin{enumerate}
\item In the $C^{*}$-setting: $V$ extends to a $C^{*}$-representation
if and only if it is isometric.
\item If $V$ is isometric, it admits a Wold decomposition $V_{n}(\zeta)=\left(S_{n}(\zeta)\tensor I_{\D}\right)\oplus Z_{n}(\zeta)$
(up to unitary equivalence) where $\D$ is some Hilbert subspace of
$\H$ and $Z$ is an isometric, fully coisometric, covariant representation
of $X$.
\item Every completely contractive, covariant representation of $X$ possesses
a minimal dilation to an isometric, covariant representation of $X$.
\end{enumerate}
\end{thm}

\begin{thm}
\label{thm:dilation_1}Let $X=\left(X(n)\right)_{n\in\Z_{+}}$ be
a standard subproduct system, and suppose that $T$ is a completely
contractive, covariant representation of $X$ on $\H$. Then there
exists a dilation $V$ of $T$ to a Hilbert space $\K$ with the following
properties:
\begin{enumerate}
\item \label{enu:dilation_1_condition_1}There exists a Hilbert space $\U$
such that $\K=(\F_{X}\tensor_{\sigma}\D)\oplus\U$.
\item \label{enu:dilation_1_condition_2}There exists a completely contractive,
fully coisometric, covariant representation $Z$ of $X$ on $\U$
such that $V_{n}(\zeta)=\left(S_{n}(\zeta)\tensor I_{\D}\right)\oplus Z_{n}(\zeta)$
for all $n\in\Z_{+}$, $\zeta\in X(n)$.
\end{enumerate}
Moreover, $\U=\left\{ 0\right\} $ if and only if $T$ is pure.\end{thm}
\begin{proof}
Let $Q:=\slim_{n\to\infty}\widetilde{T}_{n}\widetilde{T}_{n}^{*}$
(see Remark \ref{rem:T_n_tilde_T_n_tilde_adj}). From \eqref{eq:c_c_representation_condition_tilde_standard}
and \eqref{eq:c_c_representation_condition_tilde_standard_adjoints}
we obtain \begin{equation}
\widetilde{T}_{n}(I_{X(n)}\tensor Q)\widetilde{T}_{n}^{*}=\slim_{m\to\infty}\widetilde{T}_{n}(I_{X(n)}\tensor\widetilde{T}_{m}\widetilde{T}_{m}^{*})\widetilde{T}_{n}^{*}=\slim_{m\to\infty}\widetilde{T}_{n+m}\widetilde{T}_{n+m}^{*}=Q\label{eq:dilation_thm_T_n_tilde_Q}\end{equation}
(note: Remark \ref{rem:T_n_tilde_T_n_tilde_adj} implies that $I_{X(n)}\tensor Q=\slim\limits _{m\to\infty}I_{X(n)}\tensor\widetilde{T}_{m}\widetilde{T}_{m}^{*}$
because $\bigl\{\widetilde{T}_{m}\widetilde{T}_{m}^{*}\bigr\}_{m\in\Z_{+}}$
is bounded). Define $\U:=\overline{\Img Q}=\overline{\Img Q^{\frac{1}{2}}}$.
Let $Y:\H\to\U$ be the operator $Q^{\frac{1}{2}}$ with codomain
$\U$ instead of $\H$. For all $n\in\Z_{+}$, $\widetilde{T}_{n}\widetilde{T}_{n}^{*}$
belongs to $\sigma(\M)'$ by \cite[Lemma 3.6]{Tensor_algebras}. Consequently
$Q\in\sigma(\M)'$, and thus $\U$ reduces $\sigma(a)$ for all $a\in\M$.
Write $\sigma_{\U}$ for the representation of $\M$ on $\U$ sending
$a\in\M$ to $\sigma(a)_{|\U}\in B(\U)$. In the $W^{*}$-setting,
$\sigma_{\U}$ is normal. Given $n\in\Z_{+}$, form the Hilbert space
$X(n)\tensor_{\sigma_{\U}}\U$, which may be viewed as a closed subspace
of $X(n)\tensor_{\sigma}\H$ in the natural way.

Fix $n\in\Z_{+}$, and define a linear operator $\Lambda_{n}:\Img Q^{\frac{1}{2}}\to X(n)\tensor_{\sigma_{\U}}\U$
by \begin{equation}
\Lambda_{n}(Yh):=(I_{X(n)}\tensor Y)\widetilde{T}_{n}^{*}h\label{eq:dilation_thm_Lambda_n_def}\end{equation}
for all $h\in\H$. To see that this is a legitimate definition of
a contraction, note first that $Q^{\frac{1}{2}}\in\sigma(\M)'$, so
that $I_{X(n)}\tensor Y$ makes sense as a bounded operator from $X(n)\tensor_{\sigma}\H$
to $X(n)\tensor_{\sigma_{\U}}\U$; and for $h\in\H$, we obtain from
\eqref{eq:dilation_thm_T_n_tilde_Q} that\[
\bigl\Vert(I_{X(n)}\tensor Y)\widetilde{T}_{n}^{*}h\bigr\Vert^{2}=\bigl(\widetilde{T}_{n}(I_{X(n)}\tensor Q)\widetilde{T}_{n}^{*}h,h\bigr)=\bigl(Qh,h\bigr)=\left\Vert Yh\right\Vert ^{2}.\]
It follows that $\Lambda_{n}$ extends to a bounded operator from
$\U$ to $X(n)\tensor_{\sigma_{\U}}\U$, which we also denote by $\Lambda_{n}$.
From \eqref{eq:dilation_thm_Lambda_n_def} we deduce that $\Lambda_{n}Y=(I_{X(n)}\tensor Y)\widetilde{T}_{n}^{*}$.

Define $\widetilde{Z}_{n}:=\Lambda_{n}^{*}$. Then \begin{equation}
\widetilde{Z}_{n}^{*}Y=(I_{X(n)}\tensor Y)\widetilde{T}_{n}^{*}\text{ and }Y^{*}\widetilde{Z}_{n}=\widetilde{T}_{n}(I_{X(n)}\tensor Y^{*}),\label{eq:dilation_thm_Y_adj_Z_tilde}\end{equation}
thus \[
Y^{*}\widetilde{Z}_{n}\widetilde{Z}_{n}^{*}Y=\widetilde{T}_{n}(I_{X(n)}\tensor Y^{*}Y)\widetilde{T}_{n}^{*}=\widetilde{T}_{n}(I_{X(n)}\tensor Q)\widetilde{T}_{n}^{*}=Q\]
(see \eqref{eq:dilation_thm_T_n_tilde_Q}), and therefore \[
(\widetilde{Z}_{n}\widetilde{Z}_{n}^{*}Yh,Yh)=(Qh,h)=(Yh,Yh),\]
whence we conclude that $\widetilde{Z}_{n}\widetilde{Z}_{n}^{*}=I_{\U}$.
In addition, for $\zeta\in X(n)$, $a\in\M$ and $h\in\U$, we get\[
\begin{split}Y^{*}\widetilde{Z}_{n}((a\cdot\zeta)\tensor h)=\widetilde{T}_{n}((a\cdot\zeta)\tensor Y^{*}h) & =\sigma(a)\widetilde{T}_{n}(\zeta\tensor Y^{*}h)\\
 & =\sigma(a)Y^{*}\widetilde{Z}_{n}(\zeta\tensor h)=Y^{*}\sigma_{\U}(a)\widetilde{Z}_{n}(\zeta\tensor h)\end{split}
\]
(for the last equality, note that $Y\sigma(a)=\sigma_{\U}(a)Y$ due
to the commutativity of $Q$ and $\sigma(a)$; then take adjoints).
Since $Y^{*}$ is injective, we infer that $\widetilde{Z}_{n}((a\cdot\zeta)\tensor h)=\sigma_{\U}(a)\widetilde{Z}_{n}(\zeta\tensor h)$.
That is, \eqref{eq:covariance_condition_tilde} holds with $\widetilde{Z}_{n}$
in place of $\widetilde{T}$. By virtue of Lemma \ref{lem:bijection__c_c_c_rep__tilde}
there exists a completely contractive, covariant representation $(Z_{n},\sigma_{\U})$
of the $A^{*}$-correspondence $X(n)$ on $\U$ that is related to
$\widetilde{Z}_{n}$ in the usual sense of \eqref{eq:T_tilde_def}.

In a similar fashion we have (from \eqref{eq:dilation_thm_Y_adj_Z_tilde}),
for $n,m\in\Z_{+}$,\begin{multline*}
Y^{*}\widetilde{Z}_{n+m}(p_{n+m}\tensor I_{\U})_{|X(n)\tensor X(m)\tensor_{\sigma_{\U}}\U}\\
\begin{split} & =\widetilde{T}_{n+m}(I_{X(n+m)}\tensor Y^{*})(p_{n+m}\tensor I_{\U})_{|X(n)\tensor X(m)\tensor_{\sigma_{\U}}\U}\\
 & =\widetilde{T}_{n+m}(p_{n+m}\tensor I_{\H})_{|X(n)\tensor X(m)\tensor_{\sigma}\H}(I_{X(n)\tensor X(m)}\tensor Y^{*})\\
 & =\widetilde{T}_{n}(I_{X(n)}\tensor\widetilde{T}_{m})(I_{X(n)\tensor X(m)}\tensor Y^{*})\\
 & =\widetilde{T}_{n}(I_{X(n)}\tensor Y^{*}\widetilde{Z}_{m})=Y^{*}\widetilde{Z}_{n}(I_{X(n)}\tensor\widetilde{Z}_{m}),\end{split}
\end{multline*}
establishing that $\widetilde{Z}_{n+m}(p_{n+m}\tensor I_{\U})_{|X(n)\tensor X(m)\tensor_{\sigma_{\U}}\U}=\widetilde{Z}_{n}(I_{X(n)}\tensor\widetilde{Z}_{m})$
by the injectivity of $Y^{*}$. This is precisely \eqref{eq:c_c_representation_condition_tilde_standard}
with $Z$ replacing $T$. Moreover, $Y^{*}\widetilde{Z}_{0}(a\tensor h)=\widetilde{T}_{0}(a\tensor Y^{*}h)=\sigma(a)Y^{*}h=Y^{*}\sigma_{\U}(a)h$
for all $a\in\M=X(0)$ and $h\in\U$, so $Z_{0}=\sigma_{\U}$. The
requirements of Definition \ref{def:c_c_c_representation} are therefore
satisfied, making $Z=\left(Z_{n}\right)_{n\in\Z_{+}}$ a completely
contractive, covariant representation of $X$ on $\U$.

Let $\K:=(\F_{X}\tensor_{\sigma}\D)\oplus\U$. Define an operator
$W:\H\to\K$ by $W:=\left(\begin{array}{c}
K(T)\\
Y\end{array}\right)$. Then for $h\in\H$ we have from \eqref{eq:K_adj_K} \[
\left\Vert Wh\right\Vert ^{2}=\left\Vert K(T)h\right\Vert ^{2}+\left\Vert Yh\right\Vert ^{2}=\left\Vert h\right\Vert ^{2}-(Qh,h)+(Qh,h)=\left\Vert h\right\Vert ^{2},\]
that is, $W$ is an isometry. Define the sequence $V=\left(V_{n}\right)_{n\in\Z_{+}}$
of linear maps $V_{n}:X(n)\to B(\K)$ by \[
V_{n}(\zeta):=\left(\begin{array}{cc}
S_{n}(\zeta)\tensor I_{\D} & 0\\
0 & Z_{n}(\zeta)\end{array}\right)\]
for $n\in\Z_{+},\zeta\in X(n)$. Since both the induced representation
$S\tensor I_{\D}$ and $Z$ are completely contractive, covariant
representations of $X$, the family $V$ is a completely contractive,
covariant representation of $X$ on $\K$.

From \eqref{eq:dilation_thm_Y_adj_Z_tilde} one may infer by direct
calculation that $Z_{n}(\zeta)^{*}Y=YT_{n}(\zeta)^{*}$. Hence we
deduce from Proposition \ref{pro:Poisson_kernel_S_T} that \begin{equation}
V_{n}(\zeta)^{*}W=\left(\begin{array}{c}
(S_{n}(\zeta)^{*}\tensor I_{\D})K(T)\\
Z_{n}(\zeta)^{*}Y\end{array}\right)=\left(\begin{array}{c}
K(T)T_{n}(\zeta)^{*}\\
YT_{n}(\zeta)^{*}\end{array}\right)=WT_{n}(\zeta)^{*},\label{eq:dilation_thm_V_n_T_n}\end{equation}
that is, upon identifying $\H$ with its image under $W$ we have
$T_{n}(\zeta)^{*}=\left(V_{n}(\zeta)^{*}\right)_{|\H}$. The family
$V$ is, in conclusion, a dilation of $T$ to $\K$, and the proof
is complete.\end{proof}
\begin{cor}
The first part of \cite[Theorem 8.5]{Subproduct_systems_2009} follows
as a direct consequence of the last theorem, after recalling that
$\left\Vert T\right\Vert _{cb}$ in their notations equals $\bigl\Vert\widetilde{T}_{1}\bigr\Vert$
in ours (see \cite[Lemma 3.5]{Tensor_algebras}). This alternative
proof does not require the Stinespring dilation theorem.
\end{cor}

\section{Pure and relatively isometric covariant representations\label{sec:Pure-and-relatively_isometric}}

We have seen that every completely contractive, covariant representation
of a subproduct system may be {}``extended'' to a completely contractive
representation of the tensor algebra (Theorem \ref{thm:c_c_c_representation__to__A_morphism}).
Our goal in this section and the next one is to investigate when it
is possible to further extend it to a $C^{*}$-representation of the
Toeplitz algebra in the sense of the following definition.
\begin{defn}
\label{def:extend_C_representation}Let $X=\left(X(n)\right)_{n\in\Z_{+}}$
be a standard subproduct system. We shall say that a completely contractive,
covariant representation $T$ of $X$ on $\H$ \emph{extends to a
$C^{*}$-representation} if there exists a $C^{*}$-representation
$\pi$ of $\T(X)$ on $\H$ such that $\pi(S_{n}(\zeta))=T_{n}(\zeta)$
for all $n\in\Z_{+}$, $\zeta\in X(n)$.
\end{defn}
As indicated in the introduction, we split the problem into two. In
the present section, we consider chiefly \emph{pure} covariant representations.
A sufficient condition is established in Theorem \ref{thm:rel_isom_1_pure},
and it is shown that, in some important special cases, this condition
is also necessary (Proposition \ref{pro:relatively_isometric_product_system}
and Corollary \ref{cor:rel_isom_1_finite_dim_Hilbert_space}).

If $n\in\Z_{+}$, the definition of $\widetilde{S}_{n}$ implies that
$\widetilde{S}_{n}^{*}\widetilde{S}_{n}\in\L(X(n)\tensor\F_{X})$
is the projection on the subspace $X(n)\oplus X(n+1)\oplus X(n+2)\oplus\cdots$
of $X(n)\tensor\F_{X}$ (recall that $X(n+k)\subseteq X(n)\tensor X(k)$
for all $k$), and that $\widetilde{S}_{n}\widetilde{S}_{n}^{*}\in\L(\F_{X})$
is the projection on $X(n)\oplus X(n+1)\oplus X(n+2)\oplus\cdots$
(considered as a subspace of $\F_{X}$). The \emph{product} system
case (Theorem \ref{thm:isom_prod_sys}) provides further motivation
for what follows. Notice that an operator $R:\H_{1}\to\H_{2}$ is
an isometry if and only if it is a partial isometry, and $R^{*}$
is surjective.
\begin{lem}
\label{lem:relatively_isometric_preliminary}Let $X=\left(X(n)\right)_{n\in\Z_{+}}$
be a standard subproduct system, and let $T$ be a completely contractive,
covariant representation of $X$ on $\H$. If the maps $\widetilde{T}_{n}$,
$n\in\Z_{+}$, are all partial isometries, then the following conditions
are equivalent:
\begin{enumerate}
\item For all $n\in\Z_{+}$, \begin{equation}
X(n)\tensor_{\sigma}\Img\Delta_{*}\subseteq\Img\widetilde{T}_{n}^{*}.\label{eq:relatively_isometric}\end{equation}

\item For all $n\in\Z_{+}$ and $\zeta\in X(n)$, \begin{equation}
\Delta_{*}T_{n}(\zeta)^{*}T_{n}(\zeta)\Delta_{*}=\sigma(\left\langle \zeta,\zeta\right\rangle )\Delta_{*}.\label{eq:relatively_isometric_operators}\end{equation}

\end{enumerate}
\end{lem}
Observe the resemblance and difference between \eqref{eq:relatively_isometric_operators}
and \eqref{eq:intro_isometric_covariant_rep}.
\begin{proof}
Since $\widetilde{T}_{1}\widetilde{T}_{1}^{*}$ is a projection, so
is $\Delta_{*}$. Fix $n\in\Z_{+}$. Condition \eqref{eq:relatively_isometric}
holds if and only if $I_{X(n)}\tensor\Delta_{*}\leq\widetilde{T}_{n}^{*}\widetilde{T}_{n}$;
equivalently, $\bigl\Vert\widetilde{T}_{n}^{*}\widetilde{T}_{n}(\zeta\tensor\Delta_{*}h)\bigr\Vert=\bigl\Vert\zeta\tensor\Delta_{*}h\bigr\Vert$
for all $\zeta\in X(n),h\in\H$. This is exactly Condition \eqref{eq:relatively_isometric_operators},
for $\bigl\Vert\widetilde{T}_{n}^{*}\widetilde{T}_{n}(\zeta\tensor\Delta_{*}h)\bigr\Vert=\bigl\Vert\widetilde{T}_{n}(\zeta\tensor\Delta_{*}h)\bigr\Vert$
and $\widetilde{T}_{n}(\zeta\tensor\Delta_{*}h)=T_{n}(\zeta)\Delta_{*}h$
(recall that $\Delta_{*}$ commutes with the image of $\sigma$).\end{proof}
\begin{defn}
\label{def:relatively_isometric}Let $X=\left(X(n)\right)_{n\in\Z_{+}}$
be a standard subproduct system. A completely contractive, covariant
representation $T$ of $X$ is \emph{relatively isometric} if the
maps $\widetilde{T}_{n}$, $n\in\Z_{+}$, are all partial isometries,
and one (hence both) of the conditions presented in Lemma \ref{lem:relatively_isometric_preliminary}
is fulfilled.
\end{defn}
To justify this definition, we offer (apart from Theorem \ref{thm:rel_isom_1})
the following proposition, as well as Corollary \ref{cor:rel_isom_1_finite_dim_Hilbert_space}
to follow. Unless specifically told otherwise, we use the notation
$Q$ for $\slim_{n\to\infty}\widetilde{T}_{n}\widetilde{T}_{n}^{*}$
when $T$ is fixed.
\begin{prop}
\label{pro:relatively_isometric_product_system}Let $X=\left(X(n)\right)_{n\in\Z_{+}}$
be a standard \emph{product} system, and let $T$ be a pure, completely
contractive, covariant representation of $X$. Then $T$ is isometric
(in the $C^{*}$-setting this is equivalent to $T$ extending to a
$C^{*}$-representation) if and only if $T$ is relatively isometric.\end{prop}
\begin{proof}
Necessity is clear.

Sufficiency. We will use repeatedly the fact that the maps $\widetilde{T}_{n}$,
$n\in\Z_{+}$, are partial isometries, without mentioning it explicitly.
Since $T$ is pure, we need to demonstrate that $E\tensor_{\sigma}\Img(I-Q)\subseteq\Img\widetilde{T}_{1}^{*}$;
equivalently, that $E\tensor_{\sigma}\Img(\widetilde{T}_{n}\widetilde{T}_{n}^{*}-\widetilde{T}_{n+1}\widetilde{T}_{n+1}^{*})\subseteq\Img\widetilde{T}_{1}^{*}$
for all $n\in\Z_{+}$. But $\widetilde{T}_{n}\widetilde{T}_{n}^{*}-\widetilde{T}_{n+1}\widetilde{T}_{n+1}^{*}=\widetilde{T}_{n}(I_{X(n)}\tensor\Delta_{*})\widetilde{T}_{n}^{*}$
(cf.~\eqref{eq:Poisson_kernel_isometric_1}) and $\Img\widetilde{T}_{n}(I_{X(n)}\tensor\Delta_{*})\widetilde{T}_{n}^{*}\subseteq\Img\widetilde{T}_{n}(I_{X(n)}\tensor\Delta_{*})$,
so it is enough to show that \begin{equation}
\norm{\widetilde{T}_{1}\left(\zeta\tensor(\widetilde{T}_{n}(\eta\tensor\Delta_{*}h))\right)}=\norm{\zeta\tensor(\widetilde{T}_{n}(\eta\tensor\Delta_{*}h))}\label{eq:relatively_isometric_product_system_1}\end{equation}
for all $\zeta\in E$, $\eta\in E^{\tensor n}$ and $h\in\H$. The
equality $\widetilde{T}_{1}(I_{E}\tensor\widetilde{T}_{n})=\widetilde{T}_{n+1}$
holds as $X$ is a \emph{product }system, and by \eqref{eq:relatively_isometric},
the left side of \eqref{eq:relatively_isometric_product_system_1}
equals $\norm{\zeta\tensor\eta\tensor\Delta_{*}h}$, and is thus greater
or equal to the right side because $I_{E}\tensor\widetilde{T}_{n}$
is a contraction. As the reverse inequality is obvious, we are done.\end{proof}
\begin{lem}
\label{lem:finite_dim_T_tilde}Let $X=\left(X(n)\right)_{n\in\Z_{+}}$
be a standard subproduct system such that $E=X(1)$ is a Hilbert space
(thus so are all the spaces $X(n)$, $n\in\N$). Assume that $T$
is a completely contractive, covariant representation of $X$ on $\H$.
Then for $n\in\Z_{+}$ and a fixed orthonormal base $(\zeta_{i})_{I}$
of $X(n)$ we have \begin{equation}
\widetilde{T}_{n}=\left(T_{n}(\zeta_{i})\right)_{i\in I}\text{ as a row vector and }\widetilde{T}_{n}^{*}=\left(T_{n}(\zeta_{i})^{*}\right)_{i\in I}\text{ as a column vector},\label{eq:T_n_tilde_Hilbert_spaces}\end{equation}
and consequently \begin{equation}
\widetilde{T}_{n}\widetilde{T}_{n}^{*}=\sum_{i\in I}T_{n}(\zeta_{i})T_{n}(\zeta_{i})^{*},\label{eq:T_n_tilde_Hilbert_spaces_1}\end{equation}
where in \eqref{eq:T_n_tilde_Hilbert_spaces} and \eqref{eq:T_n_tilde_Hilbert_spaces_1}
the convergence is in the strong operator topology.\end{lem}
\begin{proof}
If $i,j\in I$ and $h,k\in\H$, then \[
\left\langle \zeta_{i}\tensor h,\zeta_{j}\tensor k\right\rangle =\left\langle h,\sigma(\left\langle \zeta_{i},\zeta_{j}\right\rangle )k\right\rangle =\begin{cases}
\left\langle h,k\right\rangle  & i=j\\
0 & \text{else.}\end{cases}\]
We conclude that $X(n)\tensor_{\sigma}\H=\bigoplus_{i\in I}\zeta_{i}\tensor\H$
and $\zeta_{i}\tensor\H\cong\H$ for all $i\in I$. The first part
of \eqref{eq:T_n_tilde_Hilbert_spaces} is therefore an immediate
result of the definition of $\widetilde{T}_{n}$, and the second part
follows from the first. The rest of the lemma is an easy consequence.\end{proof}
\begin{prop}
\label{pro:finite_dim_rep_implies_rel_isom}Let $X=\left(X(n)\right)_{n\in\Z_{+}}$
be a standard subproduct system such that $E=X(1)$ is a finite dimensional
Hilbert space, and suppose that $T$ is a completely contractive,
covariant representation of $X$ on $\H$. If $T$ extends to a $C^{*}$-representation
then it is relatively isometric.\end{prop}
\begin{proof}
Let $\pi:\T(X)\to B(\H)$ be a $C^{*}$-representation as in Definition
\ref{def:extend_C_representation}. We first demonstrate that $\widetilde{T}_{n}$
is a partial isometry for $n\in\Z_{+}$. This is equivalent to $\widetilde{T}_{n}\widetilde{T}_{n}^{*}\in B(\H)$
being a projection. Indeed, for a fixed orthonormal base $(\zeta_{i})_{I}$
of $X(n)$, we have from \eqref{eq:T_n_tilde_Hilbert_spaces_1} \[
\widetilde{T}_{n}\widetilde{T}_{n}^{*}=\sum_{i\in I}T_{n}(\zeta_{i})T_{n}(\zeta_{i})^{*}=\sum_{i\in I}\pi(S_{n}(\zeta_{i})S_{n}(\zeta_{i})^{*})=\pi(\sum_{i\in I}S_{n}(\zeta_{i})S_{n}(\zeta_{i})^{*})\]
(the last equality holds since all sums are \emph{finite}). But $\sum_{i\in I}S_{n}(\zeta_{i})S_{n}(\zeta_{i})^{*}=\widetilde{S}_{n}\widetilde{S}_{n}^{*}$,
which is a projection, as indicated above. Hence $\widetilde{T}_{n}\widetilde{T}_{n}^{*}$
is also a projection.

Let $(\eta_{j})_{J}$ denote an orthonormal base of $E$. From Definition
\ref{def:Delta_star} and the foregoing we deduce that $\Delta_{*}=I_{\H}-\sum_{j\in J}T_{1}(\eta_{j})T_{1}(\eta_{j})^{*}$.
Fix $n\in\Z_{+}$ and $\zeta\in X(n)$. To establish \eqref{eq:relatively_isometric_operators}
we must verify that \[
\Delta_{*}T_{n}(\zeta)^{*}T_{n}(\zeta)\Delta_{*}=\norm{\zeta}^{2}\Delta_{*}.\]
Write $\Delta_{*}^{X}:=I_{\F_{X}}-\widetilde{S}_{1}\widetilde{S}_{1}^{*}=I_{\F_{X}}-\sum_{j\in J}S_{1}(\eta_{j})S_{1}(\eta_{j})^{*}\in\T(X)$.
Since $\pi$ is a representation of $\T(X)$, one has $\Delta_{*}=\pi(\Delta_{*}^{X})$,
and it is therefore enough to ascertain that $\Delta_{*}^{X}S_{n}(\zeta)^{*}S_{n}(\zeta)\Delta_{*}^{X}=\norm{\zeta}^{2}\Delta_{*}^{X}$.
But $\Delta_{*}^{X}$ is the projection of $\F_{X}$ on the direct
summand $\C=X(0)$, so that $S_{n}(\zeta)^{*}S_{n}(\zeta)\Delta_{*}^{X}=\norm{\zeta}^{2}\Delta_{*}^{X}$,
and the proof is complete.
\end{proof}
We now present a Wold decomposition for relatively isometric covariant
representations.
\begin{thm}
\label{thm:rel_isom_1}Let $X=\left(X(n)\right)_{n\in\Z_{+}}$ be
a standard subproduct system, and suppose that $T$ is a completely
contractive, covariant representation of $X$ on $\H$. Then the following
are equivalent:
\begin{enumerate}
\item $T$ is relatively isometric.
\item There exist a Hilbert space $\U$, a unitary $W:\H\to\K:=(\F_{X}\tensor_{\sigma}\D)\oplus\U$
and a fully coisometric, covariant representation $Z=\left(Z_{n}\right)_{n\in\Z_{+}}$
of $X$ on $\U$ such that $WT_{n}(\zeta)W^{-1}=\left(S_{n}(\zeta)\tensor I_{\D}\right)\oplus Z_{n}(\zeta)$
for all $n\in\Z_{+}$, $\zeta\in X(n)$.
\end{enumerate}
Moreover, $\U$ may be chosen to be $\left\{ 0\right\} $ if and only
if $T$ is pure.\end{thm}
\begin{proof}
$(1)\Rightarrow(2)$. We examine the influence of the assumption that
$T$ is relatively isometric on several parts of the proof of Theorem
\ref{thm:dilation_1}, commencing by claiming that $K(T)$ is \emph{onto}
$\F_{X}\tensor_{\sigma}\D$. Our assumptions yield that $(\widetilde{T}_{n}\widetilde{T}_{n}^{*}-\widetilde{T}_{n+1}\widetilde{T}_{n+1}^{*})_{n\in\Z_{+}}$
is an orthogonal sequence of projections in $B(\H)$, whose sum is
$I-Q$. Since $\widetilde{T}_{1}\widetilde{T}_{1}^{*}$ is a projection,
we have $\Delta_{*}(T)=I_{\H}-\widetilde{T}_{1}\widetilde{T}_{1}^{*}$
(which is also a projection). Moreover, from \eqref{eq:Poisson_kernel_isometric_1}
we learn that $\ker\bigl((I_{X(n)}\tensor\Delta_{*}(T))\widetilde{T}_{n}^{*}\bigr)=\ker\bigl(\widetilde{T}_{n}\widetilde{T}_{n}^{*}-\widetilde{T}_{n+1}\widetilde{T}_{n+1}^{*}\bigr)$.
By \eqref{eq:relatively_isometric} we obtain $X(n)\tensor_{\sigma}\D=X(n)\tensor_{\sigma}\Img\Delta_{*}(T)\subseteq\Img\widetilde{T}_{n}^{*}$,
and consequently $X(n)\tensor_{\sigma}\D=\Img\bigl((I_{X(n)}\tensor\Delta_{*}(T))\widetilde{T}_{n}^{*}\bigr)$.
From the foregoing we conclude that $X(n)\tensor_{\sigma}\D\subseteq\Img K(T)$
for all $n\in\Z_{+}$. But $K(T)$ is a partial isometry as $K(T)^{*}K(T)=I-Q$
by \eqref{eq:K_adj_K}. Hence its range is closed, thus it equals
$\F_{X}\tensor_{\sigma}\D$.

The operator $Q$ is a projection, so that $\U=\Img Q$ and $Y:\H\to\U$
is the coisometry mapping $h\in\H$ to $Qh$. Since $(\ker K(T))^{\perp}=\Img(I-Q)=\ker Y$,
we have $\Img W=\Img K(T)\oplus\Img Y=\K$. In conclusion, $W$ is
unitary, so we can read our desired formula for $WT_{n}(\zeta)W^{-1}$
from \eqref{eq:dilation_thm_V_n_T_n}.

$(2)\Rightarrow(1)$. For $n\in\Z_{+}$, the map $\widetilde{T}_{n}$
is a partial isometry because it equals $W^{-1}\bigl((\widetilde{S}_{n}\tensor I_{\D})\oplus\widetilde{Z}_{n}\bigr)(I_{X(n)}\tensor W)$,
both $\widetilde{S}_{n}$ and $\widetilde{Z}_{n}$ are partial isometries
and $W$ is unitary. Moreover, as $Z$ is fully coisometric, we obtain
\[
\begin{split}\Delta_{*}(T) & =W^{-1}\left(\begin{array}{cc}
I_{\F_{X}\tensor_{\sigma}\D}-(\widetilde{S}_{1}\widetilde{S}_{1}^{*})\tensor I_{\D} & 0\\
0 & I_{\U}-\widetilde{Z}_{1}\widetilde{Z}_{1}^{*}\end{array}\right)W\\
 & =W^{-1}\left(\begin{array}{cc}
\Delta_{*}(S)\tensor I_{\D} & 0\\
0 & 0\end{array}\right)W.\end{split}
\]
Hence, using the fact that $\sigma(\cdot)=T_{0}(\cdot)=W^{-1}\bigl((\varphi_{\infty}(\cdot)\tensor I_{\D})\oplus Z_{0}(\cdot)\bigr)W$,
we get \[
\begin{array}{lll}
\Delta_{*}(T)T_{n}(\zeta)^{*}T_{n}(\zeta)\Delta_{*}(T) & = & W^{-1}\left(\begin{array}{cc}
\Delta_{*}(S)\tensor I_{\D} & 0\\
0 & 0\end{array}\right)\left(\begin{array}{cc}
S_{n}(\zeta)^{*}\tensor I_{\D} & 0\\
0 & Z_{n}(\zeta)^{*}\end{array}\right)\cdot\\
 &  & \cdot\left(\begin{array}{cc}
S_{n}(\zeta)\tensor I_{\D} & 0\\
0 & Z_{n}(\zeta)\end{array}\right)\left(\begin{array}{cc}
\Delta_{*}(S)\tensor I_{\D} & 0\\
0 & 0\end{array}\right)W\\
 & = & W^{-1}\left(\begin{array}{cc}
\left(\Delta_{*}(S)S_{n}(\zeta)^{*}S_{n}(\zeta)\Delta_{*}(S)\right)\tensor I_{\D} & 0\\
0 & 0\end{array}\right)W\\
 & = & W^{-1}\left(\begin{array}{cc}
\left(\varphi_{\infty}(\left\langle \zeta,\zeta\right\rangle )\Delta_{*}(S)\right)\tensor I_{\D} & 0\\
0 & 0\end{array}\right)W\\
 & = & \sigma(\left\langle \zeta,\zeta\right\rangle )\Delta_{*}(T),\end{array}\]
 proving that $T$ is relatively isometric.\end{proof}
\begin{rem*}
In the proof of $(2)\Rightarrow(1)$ we employed the operator $I_{X(n)}\tensor W$.
That it is a well-defined, bounded operator is a consequence of the
hypotheses.\end{rem*}
\begin{thm}
\label{thm:rel_isom_1_pure}Let $X=\left(X(n)\right)_{n\in\Z_{+}}$
be a standard subproduct system, and suppose that $T$ is a completely
contractive, covariant representation of $X$. If $T$ is pure and
relatively isometric, then it is induced, and it extends to a $C^{*}$-representation.\end{thm}
\begin{proof}
$T$ being pure is equivalent to being able to choose $\U=\left\{ 0\right\} $.
When this is the case, $T_{n}(\zeta)=W^{-1}(S_{n}(\zeta)\tensor I_{\D})W$
for all $n\in\Z_{+},\zeta\in X(n)$, i.e., $T$ is (unitarily equivalent
to) an induced representation. It therefore extends to a $C^{*}$-representation
in the natural fashion.\end{proof}
\begin{cor}
Let $X=\left(X(n)\right)_{n\in\Z_{+}}$ be a standard subproduct system
such that $X(n_{0})=\left\{ 0\right\} $ for some $n_{0}\in\N$. Then
every completely contractive, covariant representation of $X$ that
is relatively isometric extends to a $C^{*}$-representation.\end{cor}
\begin{proof}
Every such covariant representation is pure by Definition \ref{def:pure_ccc_representation}
because $X(n)=\left\{ 0\right\} $ for all $n\geq n_{0}$.
\end{proof}
Combining Theorem \ref{thm:rel_isom_1_pure} and Proposition \ref{pro:finite_dim_rep_implies_rel_isom}
we infer the following.
\begin{cor}
\label{cor:rel_isom_1_finite_dim_Hilbert_space}Let $X=\left(X(n)\right)_{n\in\Z_{+}}$
be a standard subproduct system such that $E=X(1)$ is a finite dimensional
Hilbert space (see Example \ref{exa:finite_dim_Hilbert_spaces}),
and suppose that $T$ is a completely contractive, covariant representation
of $X$. If $T$ is pure, then the following conditions are equivalent:
\begin{enumerate}
\item $T$ extends to a $C^{*}$-representation.
\item $T$ is relatively isometric.
\end{enumerate}
\end{cor}
We wish to emphasize the novelty of this corollary even in the well-studied
case of Arveson's symmetric subproduct system $\mathrm{SSP}_{d}$
(\cite{Arveson_subalgebras_3}; see also Example \ref{exa:fully_coisometric_symmetric_case}).
Indeed, while it has already been known that a pure, completely contractive,
covariant representation of $\mathrm{SSP}_{d}$ extends to a $C^{*}$-representation
if and only if it is an induced representation of the $\mathrm{SSP}_{d}$-shifts---the
latter is not a {}``checkable'' condition. On the other hand, our
condition of being relatively isometric is very concrete for the subproduct
systems considered in the last corollary (and particularly $\mathrm{SSP}_{d}$)
owing to Lemma \ref{lem:finite_dim_T_tilde}, or the following observation:
if $d=\dim E$ and $\left\{ e_{1},\ldots,e_{d}\right\} $ is some
orthonormal base of $E$, define $T_{i}:=T_{1}(e_{i})$, $1\leq i\leq d$.
For $\a=(\a_{1},\ldots,a_{n})\in\left\{ 1,\ldots,d\right\} ^{n}$,
set $T_{\a}:=T_{\a_{1}}\cdots T_{\a_{n}}$. Then \cite[Proposition 6.9]{Subproduct_systems_2009}
implies that \[
\widetilde{T}_{n}\widetilde{T}_{n}^{*}=\sum_{\a\in\left\{ 1,\ldots,d\right\} ^{n}}T_{\a}T_{\a}^{*}\]
for all $n\in\N$. In conclusion, the conditions in Definition \ref{def:relatively_isometric}
can be easily written in terms of the operators $T_{1},\ldots,T_{d}$.

\section{Fully coisometric covariant representations\label{sec:Fully-coisometric}}

A fully coisometric, covariant representation of a subproduct system
$X$ such that $E=X(1)$ is a finite dimensional Hilbert space possesses
a dilation to a fully coisometric, covariant representation of $X$
that extends to a $C^{*}$-representation. This is obtained by a suitable
adaptation of \cite[Theorem 2.4]{Popescu_noncomm_varieties} and \cite[Proposition 7.2]{Subproduct_systems_2009}.
The proof uses the pair Arveson's extension theorem + Stinespring's
dilation theorem, as well as properties of the Toeplitz algebra $\T(X)$,
which depend on $E$ being a finite dimensional Hilbert space. Additionally,
as Arveson's extension theorem does not supply a concrete construction
for the extension, applying the same techniques to covariant representations
of general subproduct systems, one encounters troubles when trying
to prove the $W^{*}$-setting continuity condition mentioned in Definition
\ref{def:c_c_representation_correspondence}. For these reasons, using
the same lines to prove a generalization of the aforementioned dilation
theorem seems impossible. This section is devoted to dilations and
$C^{*}$-extendability of fully coisometric, covariant representations,
complementing the developments of the previous section.

Suppose that $T$ is a completely contractive, covariant representation
of a subproduct system $X$. It is a special case of Theorem \ref{thm:c_c_c_representation__to__A_morphism}
that $\norm{T_{n}(\zeta)T_{m}(\xi)^{*}}\leq\norm{S_{n}(\zeta)S_{m}(\xi)^{*}}$
for $n,m\in\Z_{+}$, $\zeta\in X(n)$, $\xi\in X(m)$. Assume now
that $T$ extends to a $C^{*}$-representation. Since this representation
maps every {}``polynomial'' of elements of the form $S_{n}(\zeta)$
and their adjoints to the same polynomial with $T$ replacing $S$,
the norm of the latter is less that or equal to the norm of the former.
Particularly, $\norm{T_{n}(\zeta)^{*}T_{m}(\xi)}\leq\norm{S_{n}(\zeta)^{*}S_{m}(\xi)}$.
In Theorem \ref{thm:fully_coisometric_dilation} it is established
that every fully coisometric, covariant representation admits a dilation
that satisfies the last inequality. In Theorem \ref{thm:fully_coisometric_2}
we present sufficient conditions for a fully coisometric, covariant
representation to extend to a $C^{*}$-representation in the $C^{*}$-setting.
\begin{thm}
\label{thm:fully_coisometric_dilation}Let $X=\left(X(n)\right)_{n\in\Z_{+}}$
be a standard subproduct system, and suppose that $T$ is a fully
coisometric, covariant representation of $X$. Then $T$ admits a
dilation to a completely contractive, covariant representation $V$
of $X$, for which the inequality \begin{equation}
\norm{V_{n}(\zeta)^{*}V_{m}(\xi)}\leq\norm{S_{n}(\zeta)^{*}S_{m}(\xi)}\label{eq:weak_von_Neumann_inequality}\end{equation}
holds for all $n,m\in\Z_{+}$, $\zeta\in X(n)$, $\xi\in X(m)$.
\end{thm}
The motivation for using inductive limits in the following proof was
\cite[Theorem 3.7]{Quantum_Markov_Processes}, where this technique
was used for product systems with a positive \emph{real} parameter
instead of the usual Schäffer matrix construction (cf.~\cite[Theorem 2.18]{Quantum_Markov_Processes}).
\begin{proof}
We divide the proof into a few steps.

\emph{Definition of $\K$.} For $n,m\in\Z_{+}$, $n\leq m$, define
an operator $u_{n,m}:X(n)\tensor_{\sigma}\H\to X(m)\tensor_{\sigma}\H$
by \[
u_{n,m}:=(p_{m}\tensor I_{\H})(I_{X(n)}\tensor\widetilde{T}_{m-n}^{*}).\]
Then $u_{n,m}$ is a contractive $\M$-module map. Moreover, if $n\leq m\leq k$,
we have\[
\begin{split}u_{m,k}u_{n,m} & =(p_{k}\tensor I_{\H})(I_{X(m)}\tensor\widetilde{T}_{k-m}^{*})(p_{m}\tensor I_{\H})(I_{X(n)}\tensor\widetilde{T}_{m-n}^{*})\\
 & =(p_{k}\tensor I_{\H})(p_{m}\tensor I_{X(k-m)}\tensor I_{\H})(I_{X(n)}\tensor I_{X(m-n)}\tensor\widetilde{T}_{k-m}^{*})(I_{X(n)}\tensor\widetilde{T}_{m-n}^{*})\\
 & =(p_{k}\tensor I_{\H})(I_{X(n)}\tensor\widetilde{T}_{k-n}^{*})=u_{n,k}\end{split}
\]
(the passage from the second line to the third was done by virtue
of \eqref{eq:c_c_representation_condition_tilde_standard_adjoints}).
Therefore $\Bigl((X(n)\tensor_{\sigma}\H)_{n\in\Z_{+}},(u_{n,m})_{\substack{n,m\in\Z_{+}\\
n\leq m}
}\Bigr)$ is an inductive system of Hilbert spaces that are also $\M$-modules.
We define $\K$ to be its inductive limit, and denote by $u_{n}:X(n)\tensor_{\sigma}\H\to\K$
the canonical contractive $\M$-module maps satisfying $u_{m}u_{n,m}=u_{n}$
for all $n\leq m$.

Upon identifying $X(0)\tensor_{\sigma}\H=\M\tensor_{\sigma}\H$ with
$\H$ in the usual sense, we see that $u_{0,n}=\widetilde{T}_{n}^{*}$.
Since $T$ is fully coisometric, the map $u_{0,n}$ is an isometry
for each $n\in\Z_{+}$. This implies that $u_{0}:\H\to\K$ is an isometry.

In the course of the proof we shall use repeatedly the following two
analytic ({}``universality'') properties of the inductive limit:
\begin{itemize}
\item The union $\bigcup_{n\geq n_{0}}\Img u_{n}$ is dense in $\K$ for
all $n_{0}\in\Z_{+}$.
\item If $n,m\in\Z_{+}$, $x\in X(n)\tensor_{\sigma}\H$ and $y\in X(m)\tensor_{\sigma}\H$,
then \[
(u_{n}x,u_{m}y)_{\K}=\lim_{\ell\to\infty}(u_{n,\ell}x,u_{m,\ell}y)_{X(\ell)\tensor_{\sigma}\H}.\]

\end{itemize}
\emph{Definition of the operators $V_{n}(\zeta)$.} Given $n\in\Z_{+}$
and $\zeta\in X(n)$, we define an operator $V_{n}(\zeta)\in B(\K)$
by first setting \begin{equation}
V_{n}(\zeta)u_{m}(\eta\tensor h):=u_{n+m}(S_{n}(\zeta)\eta\tensor h)=u_{n+m}(p_{n+m}(\zeta\tensor\eta)\tensor h)\label{eq:fully_coisometric_dilation__V}\end{equation}
for $m\in\Z_{+}$, $\eta\in X(m)$ and $h\in\H$. We observe that
if $\ell\geq n+m$, then\begin{equation}
\begin{split}u_{n+m,\ell}(S_{n}(\zeta)\eta\tensor h) & =(p_{\ell}\tensor I_{\H})(I_{X(n+m)}\tensor\widetilde{T}_{\ell-(n+m)}^{*})(p_{n+m}(\zeta\tensor\eta)\tensor h)\\
 & =(p_{\ell}\tensor I_{\H})(p_{n+m}\tensor\widetilde{T}_{\ell-(n+m)}^{*})(\zeta\tensor\eta\tensor h)\\
 & =(p_{\ell}\tensor I_{\H})(I_{X(n)}\tensor I_{X(m)}\tensor\widetilde{T}_{\ell-(n+m)}^{*})(\zeta\tensor\eta\tensor h)\\
 & =(p_{\ell}\tensor I_{\H})(\zeta\tensor u_{m,\ell-n}(\eta\tensor h))\\
 & =(S_{n}(\zeta)\tensor I_{\H})u_{m,\ell-n}(\eta\tensor h).\end{split}
\label{eq:fully_coisometric_dilation__S_n_trick}\end{equation}
By the fact that $\norm{S_{n}(\zeta)\tensor I_{\H}}_{B(\F_{X}\tensor_{\sigma}\H)}\leq\norm{\zeta}$,
we have \begin{multline*}
\norm{u_{n+m}(S_{n}(\zeta)\eta\tensor h)}_{\K}=\lim_{\ell\to\infty}\norm{u_{n+m,\ell}(S_{n}(\zeta)\eta\tensor h)}_{X(\ell)\tensor_{\sigma}\H}\leq\\
\leq\norm{\zeta}\cdot\lim_{\ell\to\infty}\norm{u_{m,\ell-n}(\eta\tensor h)}_{X(\ell-n)\tensor_{\sigma}\H}=\norm{\zeta}\norm{u_{m}(\eta\tensor h)}_{\K}.\end{multline*}
Since the union of the ranges of the maps $u_{m}$ is dense in $\K$,
the mapping $V_{n}(\zeta)$ is well-defined, and it extends to a bounded
operator in $B(\K)$, also denoted by $V_{n}(\zeta)$.

It is important to notice that the (left) module action of $\M$ on
$\K$ is the one implemented by $V_{0}$ (to see this, simply take
$n=0$ in \eqref{eq:fully_coisometric_dilation__V}).

\emph{$V=\left(V_{n}\right)_{n\in\Z_{+}}$ is a completely contractive,
covariant representation of $X$.} The facts that $V_{0}$ is a representation
of $\M$ and that the $V_{n}$, $n\in\N$, are bimodule maps with
respect to $V_{0}$ are inferred directly from \eqref{eq:fully_coisometric_dilation__V}.
To establish that the representation $V_{0}$ is normal in the $W^{*}$-setting,
fix $m,k\in\Z_{+}$, $\eta\in X(m)$, $\xi\in X(k)$ and $h,h'\in\H$.
Then $\left(V_{0}(a)u_{m}(\eta\tensor h),u_{k}(\xi\tensor h')\right)=\left(u_{m}(a\cdot\eta\tensor h),u_{k}(\xi\tensor h')\right)$.
Since both $\sigma$ and the homomorphism that implements the left
multiplication on $X(m)$ are normal, the mapping $a\mapsto a\cdot\eta\tensor h$
is continuous when $\M$ is equipped with the ultraweak topology and
$X(m)\tensor_{\sigma}\H$ is equipped with the weak topology. On the
other hand, the linear mapping $u_{m}:X(m)\tensor_{\sigma}\H\to\K$,
being bounded, is continuous when both Hilbert spaces are equipped
with their weak topologies. Thus $a\mapsto\left(u_{m}(a\cdot\eta\tensor h),u_{k}(\xi\tensor h')\right)$
is continuous in the ultraweak topology of $\M$. This is enough to
guarantee that $V_{0}$ is normal.

If $n,m,k\in\Z_{+}$, $\zeta\in X(n)$, $\xi\in X(m)$, $\eta\in X(k)$
and $h\in\H$, then\begin{multline*}
V_{n}(\zeta)V_{m}(\xi)u_{k}(\eta\tensor h)=V_{n}(\zeta)u_{m+k}(S_{m}(\xi)\eta\tensor h)=u_{n+m+k}(S_{n}(\zeta)S_{m}(\xi)\eta\tensor h)=\\
=u_{n+m+k}(S_{n+m}(p_{n+m}(\zeta\tensor\xi))\eta\tensor h)=V_{n+m}(p_{n+m}(\zeta\tensor\xi))u_{k}(\eta\tensor h),\end{multline*}
hence $V_{n}(\zeta)V_{m}(\xi)=V_{n+m}(p_{n+m}(\zeta\tensor\xi))$.

Fix $n\in\Z_{+}$. Define an operator $\widetilde{V}_{n}:X(n)\tensor_{V_{0}}\K\to\K$
by letting \[
\widetilde{V}_{n}(\zeta\tensor u_{m}(\eta\tensor h)):=V_{n}(\zeta)u_{m}(\eta\tensor h)\]
where $\zeta,m,\eta,h$ are as usual. If $\zeta_{1},\ldots,\zeta_{p}\in X(n)$,
$m_{1},\ldots,m_{p}\in\Z_{+}$, $h_{1},\ldots h_{p}\in\H$ and $\eta_{i}\in X(m_{i})$
for $1\leq i\leq p$, then on account of \eqref{eq:fully_coisometric_dilation__V},
\eqref{eq:fully_coisometric_dilation__S_n_trick} and the fact that
$\widetilde{S}_{n}$ is contractive, we obtain\[
\begin{split}\Bigl\Vert\widetilde{V}_{n}\Bigl(\sum_{i=1}^{p}\zeta_{i}\tensor u_{m_{i}}(\eta_{i}\tensor h_{i})\Bigr)\Bigr\Vert_{\K} & =\lim_{\ell\to\infty}\Bigl\Vert\sum_{i=1}^{p}u_{n+m_{i},\ell}(S_{n}(\zeta_{i})\eta_{i}\tensor h_{i})\Bigr\Vert_{X(\ell)\tensor_{\sigma}\H}\\
 & =\lim_{\ell\to\infty}\Bigl\Vert\sum_{i=1}^{p}(S_{n}(\zeta_{i})\tensor I_{\H})u_{m_{i},\ell-n}(\eta_{i}\tensor h_{i})\Bigr\Vert_{X(\ell)\tensor_{\sigma}\H}\\
 & =\lim_{\ell\to\infty}\Bigl\Vert(\widetilde{S}_{n}\tensor I_{\H})\bigl(\sum_{i=1}^{p}\zeta_{i}\tensor u_{m_{i},\ell-n}(\eta_{i}\tensor h_{i})\bigr)\Bigr\Vert_{X(\ell)\tensor_{\sigma}\H}\\
 & \leq\lim_{\ell\to\infty}\Bigl\Vert\sum_{i=1}^{p}\zeta_{i}\tensor u_{m_{i},\ell-n}(\eta_{i}\tensor h_{i})\Bigr\Vert_{X(n)\tensor X(\ell-n)\tensor_{\sigma}\H}\\
 & =\Bigl\Vert\sum_{i=1}^{p}\zeta_{i}\tensor u_{m_{i}}(\eta_{i}\tensor h_{i})\Bigr\Vert_{X(n)\tensor_{V_{0}}\K}.\end{split}
\]
So $\widetilde{V}_{n}$ is well-defined and contractive on its domain,
and it therefore extends to a contraction with domain $X(n)\tensor_{V_{0}}\K$.
From Lemma \ref{lem:bijection__c_c_c_rep__tilde}, $V_{n}$ is a completely
contractive, covariant representation of $X(n)$ on $\K$. In particular,
the maps $V_{n}(\cdot)$ satisfy the continuity property mentioned
in Definition \ref{def:c_c_representation_correspondence} in the
$W^{*}$-setting (see Remark \ref{rem:sigma_T_normality}).

\emph{$V$ is a dilation of $T$.} Let $n\in\Z_{+}$, $\zeta\in X(n)$
be given. Then for $h\in\H$ one has \[
\begin{split}\left(u_{0}^{*}V_{n}(\zeta)u_{0}h,h\right)_{\H} & =\left(V_{n}(\zeta)u_{0}h,u_{0}h\right)_{\K}=\left(u_{n}(\zeta\tensor h),u_{0}h\right)_{\K}\\
 & =\lim_{\ell\to\infty}\left(u_{n,\ell}(\zeta\tensor h),u_{0,\ell}h\right)_{X(\ell)\tensor_{\sigma}\H}\\
 & =\lim_{\ell\to\infty}\bigl((p_{\ell}\tensor I_{\H})(I_{X(n)}\tensor\widetilde{T}_{\ell-n}^{*})(\zeta\tensor h),\widetilde{T}_{\ell}^{*}h\bigr)_{X(\ell)\tensor_{\sigma}\H}\\
 & =\lim_{\ell\to\infty}\bigl(\widetilde{T}_{\ell}(p_{\ell}\tensor I_{\H})(I_{X(n)}\tensor\widetilde{T}_{\ell-n}^{*})(\zeta\tensor h),h\bigr)_{\H}.\end{split}
\]
 But $\widetilde{T}_{\ell}(p_{\ell}\tensor I_{\H})=\widetilde{T}_{n}(I_{X(n)}\tensor\widetilde{T}_{\ell-n})$
(by \eqref{eq:c_c_representation_condition_tilde_standard}) and $T$
is fully coisometric, so we conclude that $\left(u_{0}^{*}V_{n}(\zeta)u_{0}h,h\right)=\left(T_{n}(\zeta)h,h\right)$.
Consequently, $u_{0}^{*}V_{n}(\zeta)u_{0}=T_{n}(\zeta)$.

If, moreover, $m\in\Z_{+}$, $\eta\in X(m)$ and $h\in\H$, then by
repeated usage of the foregoing, \[
\begin{split}u_{0}^{*}V_{n}(\zeta)u_{m}(\eta\tensor h) & =u_{0}^{*}V_{n}(\zeta)V_{m}(\eta)u_{0}h=u_{0}^{*}V_{n+m}(p_{n+m}(\zeta\tensor\eta))u_{0}h\\
 & =T_{n+m}(p_{n+m}(\zeta\tensor\eta))h=T_{n}(\zeta)T_{m}(\eta)h\\
 & =u_{0}^{*}V_{n}(\zeta)u_{0}u_{0}^{*}V_{m}(\eta)u_{0}h=u_{0}^{*}V_{n}(\zeta)u_{0}u_{0}^{*}u_{m}(\eta\tensor h).\end{split}
\]
 Therefore $u_{0}^{*}V_{n}(\zeta)=u_{0}^{*}V_{n}(\zeta)u_{0}u_{0}^{*}$,
which implies that $(u_{0}\H)^{\perp}$ is invariant under $V_{n}(\zeta)$
because $u_{0}u_{0}^{*}$ is a projection. In conclusion, $V$ is
a dilation of $T$.

\emph{Proof of \eqref{eq:weak_von_Neumann_inequality}.} Fix \emph{$n,m\in\Z_{+}$,
$\zeta\in X(n)$, $\xi\in X(m)$.} Given $p,q,t\in\Z_{+}$, $\mu_{1},\ldots,\mu_{t}\in X(p)$,
$\nu_{1},\ldots,\nu_{t}\in X(q)$ and $h_{1},\ldots,h_{t},h_{1}',\ldots h_{t}'\in\H$,
denote $x:=u_{p}(\sum_{i=1}^{t}\mu_{i}\tensor h_{i})$, $y:=u_{q}(\sum_{j=1}^{t}\nu_{j}\tensor h_{j}')$.
Then from \eqref{eq:fully_coisometric_dilation__V} and \eqref{eq:fully_coisometric_dilation__S_n_trick}
we obtain \begin{multline*}
\left(V_{n}(\zeta)^{*}V_{m}(\xi)x,y\right)=\left(V_{m}(\xi)x,V_{n}(\zeta)y\right)=\\
\begin{split} & =\sum_{i,j=1}^{t}\left(u_{p+m}(S_{m}(\xi)\mu_{i}\tensor h_{i}),u_{q+n}(S_{n}(\zeta)\nu_{j}\tensor h_{j}')\right)\\
 & =\lim_{\ell\to\infty}\sum_{i,j=1}^{t}\left((S_{m}(\xi)\tensor I_{\H})u_{p,\ell-m}(\mu_{i}\tensor h_{i}),(S_{n}(\zeta)\tensor I_{\H})u_{q,\ell-n}(\nu_{j}\tensor h_{j}')\right)\\
 & =\lim_{\ell\to\infty}\Bigl((S_{n}(\zeta)^{*}S_{m}(\xi)\tensor I_{\H})\bigl(\sum_{i=1}^{t}u_{p,\ell-m}(\mu_{i}\tensor h_{i})\bigr),\bigl(\sum_{j=1}^{t}u_{q,\ell-n}(\nu_{j}\tensor h_{j}')\bigr)\Bigr).\end{split}
\end{multline*}
Hence \begin{multline*}
\left|\left(V_{n}(\zeta)^{*}V_{m}(\xi)x,y\right)\right|\leq\\
\begin{split} & \leq\norm{S_{n}(\zeta)^{*}S_{m}(\xi)\tensor I_{\H}}\cdot\lim_{\ell\to\infty}\Bigl\Vert u_{p,\ell-m}(\sum_{i=1}^{t}\mu_{i}\tensor h_{i})\Bigr\Vert\cdot\lim_{\ell\to\infty}\Bigl\Vert u_{q,\ell-n}(\sum_{j=1}^{t}\nu_{j}\tensor h_{j}')\Bigr\Vert\\
 & =\norm{S_{n}(\zeta)^{*}S_{m}(\xi)\tensor I_{\H}}\cdot\norm x\cdot\norm y\leq\norm{S_{n}(\zeta)^{*}S_{m}(\xi)}\cdot\norm x\cdot\norm y.\end{split}
\end{multline*}
Since $x,y$ are arbitrary elements of a dense subset of $\K$, inequality
\eqref{eq:weak_von_Neumann_inequality} follows\emph{.}\end{proof}
\begin{rem}
We do not know whether the dilation $V$ constructed in the theorem
is automatically fully coisometric. Nevertheless, we do know from
\eqref{eq:fully_coisometric_dilation__V} that $\widetilde{V}_{n}$
has dense range for all $n$.\end{rem}
\begin{thm}
\label{thm:fully_coisometric_2}Let $X=\left(X(n)\right)_{n\in\Z_{+}}$
be a standard subproduct system in the $C^{*}$-setting, and suppose
that $T$ is a fully coisometric, covariant representation of $X$
on $\H$, satisfying \begin{equation}
\lim_{\ell\to\infty}\bigl\Vert(p_{\ell}\tensor I_{\H})(\eta\tensor\widetilde{T}_{\ell-m}^{*}h)\bigr\Vert_{X(\ell)\tensor_{\sigma}\H}=\norm{T_{m}(\eta)h}_{\H}\label{eq:fully_coisometric_C_star_representation}\end{equation}
for all $m\in\N$, $\eta\in X(m)$ and $h\in\H$. Then $T$ extends
to a $C^{*}$-representation.\end{thm}
\begin{rem}
If $m,\eta,h$ are as above and $\ell\geq m$, then

\begin{multline*}
\bigl\Vert(p_{\ell+1}\tensor I_{\H})(\eta\tensor\widetilde{T}_{\ell+1-m}^{*}h)\bigr\Vert_{X(\ell+1)\tensor_{\sigma}\H}^{2}=\\
\begin{array}{lll}
 & = & \bigl((I_{X(m)}\tensor\widetilde{T}_{\ell+1-m})(p_{\ell+1}\tensor I_{\H})(I_{X(m)}\tensor\widetilde{T}_{\ell+1-m}^{*})(\eta\tensor h),\eta\tensor h\bigr)\\
 & = & \bigl((I_{X(m)}\tensor\widetilde{T}_{\ell-m})(I_{X(\ell)}\tensor\widetilde{T}_{1})(p_{\ell+1}\tensor I_{\H})(I_{X(m)}\tensor I_{X(\ell-m)}\tensor\widetilde{T}_{1}^{*})\\
 &  & (I_{X(m)}\tensor\widetilde{T}_{\ell-m}^{*})(\eta\tensor h),\eta\tensor h\bigr)\\
 & = & \bigl((I_{X(m)}\tensor\widetilde{T}_{\ell-m})(I_{X(\ell)}\tensor\widetilde{T}_{1})(p_{\ell+1}\tensor I_{\H})(I_{X(\ell)}\tensor\widetilde{T}_{1}^{*})(p_{\ell}\tensor I_{\H})\\
 &  & (I_{X(m)}\tensor\widetilde{T}_{\ell-m}^{*})(\eta\tensor h),\eta\tensor h\bigr)\\
 & \leq & \bigl((I_{X(m)}\tensor\widetilde{T}_{\ell-m})(p_{\ell}\tensor I_{\H})(I_{X(m)}\tensor\widetilde{T}_{\ell-m}^{*})(\eta\tensor h),\eta\tensor h\bigr)\\
 & = & \bigl\Vert(p_{\ell}\tensor I_{\H})(\eta\tensor\widetilde{T}_{\ell-m}^{*}h)\bigr\Vert_{X(\ell)\tensor_{\sigma}\H}^{2}.\end{array}\end{multline*}
That is, the sequence $\bigl\{\bigl\Vert(p_{\ell}\tensor I_{\H})(\eta\tensor\widetilde{T}_{\ell-m}^{*}h)\bigr\Vert\bigr\}_{\ell\geq m}$
is decreasing. Furthermore, \[
\begin{split}\bigl\Vert T_{m}(\eta)h\bigr\Vert^{2} & =\bigl(\widetilde{T}_{m}^{*}\widetilde{T}_{m}(\eta\tensor h),\eta\tensor h\bigr)\\
 & =\bigl((I_{X(m)}\tensor\widetilde{T}_{\ell-m})\widetilde{T}_{\ell}^{*}\widetilde{T}_{\ell}(p_{\ell}\tensor I_{\H})(I_{X(m)}\tensor\widetilde{T}_{\ell-m}^{*})(\eta\tensor h),\eta\tensor h\bigr),\end{split}
\]
whence (since $I_{X(\ell)\tensor_{\sigma}\H}-\widetilde{T}_{\ell}^{*}\widetilde{T}_{\ell}$
is a projection) \begin{multline*}
\bigl\Vert(p_{\ell}\tensor I_{\H})(\eta\tensor\widetilde{T}_{\ell-m}^{*}h)\bigr\Vert^{2}-\bigl\Vert T_{m}(\eta)h\bigr\Vert^{2}=\\
\begin{split} & =\bigl((I_{X(m)}\tensor\widetilde{T}_{\ell-m})(I_{X(\ell)\tensor_{\sigma}\H}-\widetilde{T}_{\ell}^{*}\widetilde{T}_{\ell})(p_{\ell}\tensor I_{\H})(I_{X(m)}\tensor\widetilde{T}_{\ell-m}^{*})(\eta\tensor h),\eta\tensor h\bigr)\\
 & =\bigl\Vert(I_{X(\ell)\tensor_{\sigma}\H}-\widetilde{T}_{\ell}^{*}\widetilde{T}_{\ell})(p_{\ell}\tensor I_{\H})(I_{X(m)}\tensor\widetilde{T}_{\ell-m}^{*})(\eta\tensor h)\bigr\Vert^{2}\\
 & =\bigl\Vert(p_{\ell}\tensor I_{\H})(\eta\tensor\widetilde{T}_{\ell-m}^{*}h)-\widetilde{T}_{\ell}^{*}T_{m}(\eta)h\bigr\Vert^{2}\geq0.\end{split}
\end{multline*}
In conclusion, the limit on the left side of \eqref{eq:fully_coisometric_C_star_representation}
always exists, and is greater than or equal to the right side. Moreover,
\eqref{eq:fully_coisometric_C_star_representation} is equivalent
to \begin{equation}
\lim_{\ell\to\infty}\bigl\Vert(p_{\ell}\tensor I_{\H})(\eta\tensor\widetilde{T}_{\ell-m}^{*}h)-\widetilde{T}_{\ell}^{*}T_{m}(\eta)h\bigr\Vert_{X(\ell)\tensor_{\sigma}\H}=0.\label{eq:fully_coisometric_C_star_representation_2}\end{equation}

\end{rem}
Before proving the theorem we turn our attention to condition \eqref{eq:fully_coisometric_C_star_representation},
showing that in the case of \emph{product} systems, it is equivalent
to $T$ being isometric. This is the analogue of Proposition \ref{pro:relatively_isometric_product_system}
to the fully coisometric case. The condition is further discussed
in Example \ref{exa:fully_coisometric_symmetric_case}.
\begin{prop}
Let $X=\left(X(n)\right)_{n\in\Z_{+}}$ be a standard \emph{product}
system, and suppose that $T$ is a fully coisometric, covariant representation
of $X$. Then $T$ is isometric (in the $C^{*}$-setting this is equivalent
to $T$ extending to a $C^{*}$-representation) if and only if $T$
satisfies \eqref{eq:fully_coisometric_C_star_representation}.\end{prop}
\begin{proof}
Since $p_{\ell}=I_{E^{\tensor\ell}}$ and the operators $I_{X(m)}\tensor\widetilde{T}_{\ell-m}^{*}$
are isometric for every $\ell\geq m$, \eqref{eq:fully_coisometric_C_star_representation}
is equivalent to the equality $\left\Vert \eta\tensor h\right\Vert =\bigl\Vert\widetilde{T}_{m}(\eta\tensor h)\bigr\Vert$.
This is true for all $m,\eta,h$ if and only if $T$ is isometric.
\end{proof}
We now return to Theorem \ref{thm:fully_coisometric_2} and its proof.
\begin{defn}
Let $X=\left(X(n)\right)_{n\in\Z_{+}}$ be a standard subproduct system.
An $S$\emph{-monomial} is a composition of finitely many of the operators
$S_{n}(\zeta)$ ($n\in\Z_{+}$, $\zeta\in X(n)$) and their adjoints.
Every such operator can be written as $\prod_{i=1}^{t}S_{m_{i}}(\xi_{i})^{*}S_{n_{i}}(\zeta_{i})$
for suitable $t\in\N$ and $n_{i},m_{i}\in\Z_{+}$, $\zeta_{i}\in X(n_{i})$,
$\xi_{i}\in X(m_{i})$ $(1\leq i\leq t)$ (in case $\M$ is not unital,
it may be implicitly replaced by $\M^{1}$ where necessary, letting
$S_{0}(I):=I_{\F_{X}}$). An $S$-monomial of this form is\emph{ }said
to be of\emph{ degree }$\sum_{i=1}^{t}(n_{i}-m_{i})$. For $k\in\Z$,
define $\T_{k}(X)$ to be the closed linear span of all $S$-monomials
of degree $k$. Evidently, $\T_{k}(X)$ is an operator space contained
in $\T(X)$. In the special case $k=0$, $\T_{0}(X)$ is a $C^{*}$-subalgebra
of $\T(X)$.\end{defn}
\begin{lem}
\label{lem:fully_coisometric_2}Under the assumptions of Theorem \ref{thm:fully_coisometric_2}
there exists, to each $k\in\Z$, a contraction $\pi_{k}:\T_{k}(X)\to B(\H)$
mapping an $S$-monomial $\prod_{i=1}^{t}S_{m_{i}}(\xi_{i})^{*}S_{n_{i}}(\zeta_{i})$
of degree $k$ to $\prod_{i=1}^{t}T_{m_{i}}(\xi_{i})^{*}T_{n_{i}}(\zeta_{i})$.
Moreover, $\pi_{0}$ is a $C^{*}$-representation of $\T_{0}(X)$
on $\H$.\end{lem}
\begin{proof}
Returning to the proof of Theorem \ref{thm:fully_coisometric_dilation},
we claim that \eqref{eq:fully_coisometric_C_star_representation}
yields that $\K=u_{0}\H$. Indeed, let $m\in\N$, $\eta\in X(m)$
and $h\in\H$ be given, and write $g:=T_{m}(\eta)h\in\H$. Then from
\eqref{eq:fully_coisometric_C_star_representation_2} we deduce that
\[
0=\lim_{\ell\to\infty}\norm{u_{m,\ell}(\eta\tensor h)-u_{0,\ell}g}=\lim_{\ell\to\infty}\norm{u_{m,\ell}(\eta\tensor h-u_{0,m}g)}.\]
As a result, $u_{m}(\eta\tensor h)=u_{m}u_{0,m}g=u_{0}g\in u_{0}\H$,
and we conclude that $\K=u_{0}\H$. Since $u_{0}$ is isometric it
is, in fact, unitary. In Theorem \ref{thm:fully_coisometric_dilation}
we have established the existence of a dilation $V$ of $T$, satisfying
\begin{equation}
u_{0}^{*}V_{n}(\zeta)u_{0}=T_{n}(\zeta)\label{eq:fully_coisometric_2__T_V}\end{equation}
for $n\in\Z_{+}$ and $\zeta\in X(n)$. Consequently, under current
circumstances, the operator $V_{n}(\zeta)$ is unitarily equivalent
to $T_{n}(\zeta)$.

If $n\in\Z_{+}$, $\zeta\in X(n)$ and $h\in\H$, then $V_{n}(\zeta)u_{0}h=u_{n}(\zeta\tensor h)$,
and from \eqref{eq:fully_coisometric_dilation__S_n_trick} (or simply
the definition of $u_{n,\ell}$) we have $u_{n,\ell}(\zeta\tensor h)=(S_{n}(\zeta)\tensor I_{\H})u_{0,\ell-n}h$.
On the other hand, $u_{n}(\zeta\tensor h)=u_{0}x$ for some $x\in\H$,
hence $u_{n}(\zeta\tensor h-u_{0,n}x)=0$. This is equivalent to $\lim_{\ell\to\infty}\norm{u_{n,\ell}(\zeta\tensor h-u_{0,n}x)}=0$,
i.e.,\begin{equation}
V_{n}(\zeta)u_{0}h=u_{0}x\text{ and }\lim_{\ell\to\infty}\norm{u_{0,\ell}x-(S_{n}(\zeta)\tensor I_{\H})u_{0,\ell-n}h}=0.\label{eq:fully_coisometric_2__V_S}\end{equation}

We now demonstrate a similar relation between $V_{n}(\zeta)^{*}$
and $S_{n}(\zeta)^{*}\tensor I_{\H}$. Given $n\in\Z_{+}$, $\zeta\in X(n)$
and $h\in\H$, consider $y:=(S_{n}(\zeta)^{*}\tensor I_{\H})u_{0,n}h=(S_{n}(\zeta)^{*}\tensor I_{\H})\widetilde{T}_{n}^{*}h\in\M\tensor\H$.
If $\ell\geq0$, then \[
\begin{split}u_{0,\ell}y & =(p_{\ell}\tensor I_{\H})(I_{\M}\tensor\widetilde{T}_{\ell}^{*})(S_{n}(\zeta)^{*}\tensor I_{\H})\widetilde{T}_{n}^{*}h\\
 & =(p_{\ell}\tensor I_{\H})(S_{n}(\zeta)^{*}\tensor I_{X(\ell)\tensor_{\sigma}\H})(I_{X(n)}\tensor\widetilde{T}_{\ell}^{*})\widetilde{T}_{n}^{*}h\\
 & \overset{(\#)}{=}(S_{n}(\zeta)^{*}\tensor I_{\H})\widetilde{T}_{n+\ell}^{*}h=(S_{n}(\zeta)^{*}\tensor I_{\H})u_{0,n+\ell}h.\end{split}
\]
Equality (\#) holds because the vector on which $p_{\ell}\tensor I_{\H}$
acts already belongs to ($\M\tensor X(\ell)\tensor_{\sigma}\H$, which
we identify with) $X(\ell)\tensor_{\sigma}\H$. This implies that
\[
\begin{split}(V_{n}(\zeta)^{*}u_{0}h,u_{0}g)_{\K} & =(u_{0}h,V_{n}(\zeta)u_{0}g)_{\K}\\
 & =\lim_{\ell\to\infty}(u_{0,n+\ell}h,(S_{n}(\zeta)\tensor I_{\H})u_{0,\ell}g)_{X(n+\ell)\tensor_{\sigma}\H}\\
 & =\lim_{\ell\to\infty}((S_{n}(\zeta)^{*}\tensor I_{\H})u_{0,n+\ell}h,u_{0,\ell}g)_{X(\ell)\tensor_{\sigma}\H}\\
 & =\lim_{\ell\to\infty}(u_{0,\ell}y,u_{0,\ell}g)_{X(\ell)\tensor_{\sigma}\H}=(u_{0}y,u_{0}g)_{\K}.\end{split}
\]
In conclusion, \begin{equation}
V_{n}(\zeta)^{*}u_{0}h=u_{0}y\text{ and }u_{0,\ell}y=(S_{n}(\zeta)^{*}\tensor I_{\H})u_{0,\ell+n}h\text{ for all }\ell.\label{eq:fully_coisometric_2__V_S_stars}\end{equation}

Let us prove that \begin{equation}
\Bigl\Vert\prod_{i=1}^{t}T_{m_{i}}(\xi_{i})^{*}T_{n_{i}}(\zeta_{i})\Bigr\Vert\leq\Bigl\Vert\prod_{i=1}^{t}S_{m_{i}}(\xi_{i})^{*}S_{n_{i}}(\zeta_{i})\Bigr\Vert\label{eq:fully_coisometric_2}\end{equation}
for every $t\in\N$ and $n_{i},m_{i}\in\Z_{+}$, $\zeta_{i}\in X(n_{i})$,
$\xi_{i}\in X(m_{i})$ $(1\leq i\leq t)$. Till the end of the proof,
the symbol $\prod$ will stand for multiplication in reverse order,
for the sake of convenience. From \eqref{eq:fully_coisometric_2__T_V}
and $u_{0}$ being unitary, \eqref{eq:fully_coisometric_2} is equivalent
to the inequality \[
\Bigl\Vert\prod_{i=1}^{t}V_{m_{i}}(\xi_{i})^{*}V_{n_{i}}(\zeta_{i})u_{0}\Bigr\Vert\leq\Bigl\Vert\prod_{i=1}^{t}S_{m_{i}}(\xi_{i})^{*}S_{n_{i}}(\zeta_{i})\Bigr\Vert.\]
Fix $t$ and $n_{i},m_{i},\zeta_{i},\xi_{i}$ ($1\leq i\leq t$) as
needed, and let $h\in\H$. We prove inductively that for all $0\leq p\leq t$
there exists $y_{p}\in\H$ such that \begin{equation}
\begin{split} & \quad\prod_{i=1}^{p}V_{m_{i}}(\xi_{i})^{*}V_{n_{i}}(\zeta_{i})u_{0}h=u_{0}y_{p}\text{ and }\\
 & \lim_{\ell\to\infty}\Bigl\Vert u_{0,\ell}y_{p}-\Bigl(\prod_{i=1}^{p}S_{m_{i}}(\xi_{i})^{*}S_{n_{i}}(\zeta_{i})\tensor I_{\H}\Bigr)u_{0,\ell-\sum_{i=1}^{p}(n_{i}-m_{i})}h\Bigr\Vert=0.\end{split}
\label{eq:fully_coisometric_2__induction_assumption}\end{equation}
For $p=0$ choose $y_{0}:=h$, for which \eqref{eq:fully_coisometric_2__induction_assumption}
surely holds. Assuming the existence of $y_{p}$ has been exhibited,
let $x_{p+1}\in\H$ be the unique vector such that $u_{0}x_{p+1}=V_{n_{p+1}}(\zeta_{p+1})\prod_{i=1}^{p}V_{m_{i}}(\xi_{i})^{*}V_{n_{i}}(\zeta_{i})u_{0}h=V_{n_{p+1}}(\zeta_{p+1})u_{0}y_{p}$.
Then \eqref{eq:fully_coisometric_2__V_S} implies that \[
\lim_{\ell\to\infty}\norm{u_{0,\ell}x_{p+1}-(S_{n_{p+1}}(\zeta_{p+1})\tensor I_{\H})u_{0,\ell-n_{p+1}}y_{p}}=0.\]
Consequently, from \eqref{eq:fully_coisometric_2__induction_assumption}
and the boundedness of $S_{n_{p+1}}(\zeta_{p+1})\tensor I_{\H}$,
\begin{equation}
\lim_{\ell\to\infty}\Bigl\Vert u_{0,\ell}x_{p+1}-\Bigl(S_{n_{p+1}}(\zeta_{p+1})\prod_{i=1}^{p}S_{m_{i}}(\xi_{i})^{*}S_{n_{i}}(\zeta_{i})\tensor I_{\H}\Bigr)u_{0,\ell-n_{p+1}-\sum_{i=1}^{p}(n_{i}-m_{i})}h\Bigr\Vert=0.\label{eq:fully_coisometric_2__induction_step1}\end{equation}
Denote by $y_{p+1}$ the element of $\H$ satisfying $u_{0}y_{p+1}=\prod_{i=1}^{p+1}V_{m_{i}}(\xi_{i})^{*}V_{n_{i}}(\zeta_{i})u_{0}h=V_{m_{p+1}}(\xi_{p+1})^{*}u_{0}x_{p+1}$.
By \eqref{eq:fully_coisometric_2__V_S_stars} we deduce that \[
u_{0,\ell}y_{p+1}=(S_{m_{p+1}}(\xi_{p+1})^{*}\tensor I_{\H})u_{0,\ell+m_{p+1}}x_{p+1}.\]
The operator $S_{m_{p+1}}(\xi_{p+1})^{*}\tensor I_{\H}$ being bounded,
we conclude in light of \eqref{eq:fully_coisometric_2__induction_step1}
that \[
\lim_{\ell\to\infty}\Bigl\Vert u_{0,\ell}y_{p+1}-\Bigl(\prod_{i=1}^{p+1}S_{m_{i}}(\xi_{i})^{*}S_{n_{i}}(\zeta_{i})\tensor I_{\H}\Bigr)u_{0,\ell-\sum_{i=1}^{p+1}(n_{i}-m_{i})}h\Bigr\Vert=0.\]
So by induction, \eqref{eq:fully_coisometric_2__induction_assumption}
holds with $p=t$. Therefore \begin{multline}
\Bigl\Vert\prod_{i=1}^{t}V_{m_{i}}(\xi_{i})^{*}V_{n_{i}}(\zeta_{i})u_{0}h\Bigr\Vert=\norm{u_{0}y_{t}}=\lim_{\ell\to\infty}\norm{u_{0,\ell}y_{t}}\\
\begin{split} & =\lim_{\ell\to\infty}\Bigl\Vert\Bigl(\prod_{i=1}^{t}S_{m_{i}}(\xi_{i})^{*}S_{n_{i}}(\zeta_{i})\tensor I_{\H}\Bigr)u_{0,\ell-\sum_{i=1}^{t}(n_{i}-m_{i})}h\Bigr\Vert\\
 & \leq\Bigl\Vert\prod_{i=1}^{t}S_{m_{i}}(\xi_{i})^{*}S_{n_{i}}(\zeta_{i})\tensor I_{\H}\Bigr\Vert\cdot\lim_{\ell\to\infty}\norm{u_{0,\ell-\sum_{i=1}^{t}(n_{i}-m_{i})}h}\\
 & \leq\Bigl\Vert\prod_{i=1}^{t}S_{m_{i}}(\xi_{i})^{*}S_{n_{i}}(\zeta_{i})\Bigr\Vert\cdot\norm h,\end{split}
\label{eq:fully_coisometric_2__V_u_0_S}\end{multline}
proving \eqref{eq:fully_coisometric_2}. If $k\in\Z$, the existence
of the mapping $\pi_{k}$ is shown similarly. Assume that $\prod_{i=1}^{t}S_{m_{i}^{(j)}}(\xi_{i}^{(j)})^{*}S_{n_{i}^{(j)}}(\zeta_{i}^{(j)})$
($j=1,\ldots,q$) is a finite family of $S$-monomials of degree $k$,
and fix $h\in\H$. Using \eqref{eq:fully_coisometric_2__induction_assumption}
$q$ times, we furnish the existence of $y\in\H$ such that $\sum_{j=1}^{q}\prod_{i=1}^{t}V_{m_{i}^{(j)}}(\xi_{i}^{(j)})^{*}V_{n_{i}^{(j)}}(\zeta_{i}^{(j)})u_{0}h=u_{0}y$
and \[
\lim_{\ell\to\infty}\Bigl\Vert u_{0,\ell}y-\Bigl(\sum_{j=1}^{q}\prod_{i=1}^{t}S_{m_{i}^{(j)}}(\xi_{i}^{(j)})^{*}S_{n_{i}^{(j)}}(\zeta_{i}^{(j)})\tensor I_{\H}\Bigr)u_{0,\ell-k}h\Bigr\Vert=0.\]
Consequently, just as in \eqref{eq:fully_coisometric_2__V_u_0_S},
we have

\begin{multline*}
\Bigl\Vert\sum_{j=1}^{q}\prod_{i=1}^{t}V_{m_{i}^{(j)}}(\xi_{i}^{(j)})^{*}V_{n_{i}^{(j)}}(\zeta_{i}^{(j)})u_{0}h\Bigr\Vert\leq\\
\begin{split} & \leq\Bigl\Vert\sum_{j=1}^{q}\prod_{i=1}^{t}S_{m_{i}^{(j)}}(\xi_{i}^{(j)})^{*}S_{n_{i}^{(j)}}(\zeta_{i}^{(j)})\tensor I_{\H}\Bigr\Vert\cdot\lim_{\ell\to\infty}\norm{u_{0,\ell-k}h}\\
 & \leq\Bigl\Vert\sum_{j=1}^{q}\prod_{i=1}^{t}S_{m_{i}^{(j)}}(\xi_{i}^{(j)})^{*}S_{n_{i}^{(j)}}(\zeta_{i}^{(j)})\Bigr\Vert\cdot\norm h.\end{split}
\end{multline*}
In conclusion, the canonical linear mapping from the span of all $S$-monomials
of degree $k$ to $B(\H)$, mapping $\prod_{i=1}^{t}S_{m_{i}}(\xi_{i})^{*}S_{n_{i}}(\zeta_{i})$
to $\prod_{i=1}^{t}T_{m_{i}}(\xi_{i})^{*}T_{n_{i}}(\zeta_{i})$, is
a well-defined contraction, which is a multiplicative $*$-mapping
when $k=0$. It therefore extends to a contraction $\pi_{k}$ from
$\T_{k}(X)$ to $B(\H)$. Since $\T_{0}(X)$ is a $C^{*}$-algebra,
$\pi_{0}$ is a $C^{*}$-representation.\end{proof}
\begin{rem}
Condition \eqref{eq:fully_coisometric_C_star_representation} is not
only sufficient, but also necessary, to having $\K=u_{0}\H$. We omit
the details.\end{rem}
\begin{proof}
[Proof of Theorem \ref{thm:fully_coisometric_2}]Let us begin with
some facts on circle actions. Suppose that $B$ is a $C^{*}$-algebra
with an action $\a$ of $\mathbb{T}$ on $B$. Denote the spectral
subspaces for $\a$ (\cite[Definition 2.1]{Exel_Circle_Actions})
by $\left(B_{k}\right)_{k\in\Z}$. Then $B_{1}$ becomes a $B_{0}-B_{0}$
Hilbert $C^{*}$-bimodule in the sense of \cite[Definition 1.8]{Brown_Mingo_Shen}
upon letting $\left\langle a,b\right\rangle _{R}:=a^{*}b$ and $\left\langle a,b\right\rangle _{L}:=ab^{*}$
for all $a,b\in B_{1}$. The norm in $B_{1}$ as a bimodule is the
same as its natural one. If now $\a$ is semi-saturated, i.e., $B$
is generated as a $C^{*}$-algebra by $B_{0}$ and $B_{1}$ (see \cite[Definition 4.1]{Exel_Circle_Actions}),
then by \cite[Theorem 3.1]{Abadie_Eilers_Exel} and its proof we have
$B\cong B_{0}\rtimes_{B_{1}}\Z$, the crossed product being the one
defined in \cite[Definition 2.4]{Abadie_Eilers_Exel}; and moreover,
if the crossed product maps are denoted by $\iota_{B_{0}},\iota_{B_{1}}$,
then the implementing isomorphism, say $\phi$, makes the following
diagram commute:\[
\xymatrix{(B_{0},B_{1})\ar[d]_{(\iota_{B_{0}},\iota_{B_{1}})}\ar[rd]\\
B_{0}\rtimes_{B_{1}}\Z\ar[r]\sb-\phi & B}
\]
(the maps on the diagonal being the inclusions).

In our framework, we construct the usual gauge action of $\mathbb{T}$
on $\T(X)$. Given $\lambda\in\mathbb{T}$, define a unitary $W_{\lambda}\in\L(\F_{X})$
by $\bigoplus_{n\in\Z_{+}}\zeta_{n}\mapsto\bigoplus_{n\in\Z_{+}}\lambda^{n}\zeta_{n}$.
From the definition of the $X$-shift it follows that if $A$ is an
$S$-monomial of degree $k\in\Z$ in $\T(X)$, then \begin{equation}
W_{\lambda}AW_{\l}^{*}=\l^{k}A.\label{eq:fully_coisometric_2__gauge}\end{equation}
Hence, the formula $\a_{\l}(A):=W_{\lambda}AW_{\l}^{*}$ defines an
automorphism of $\T(X)$, and $\l\mapsto\a_{\l}$ is an action of
$\mathbb{T}$ on $\T(X)$ (in particular, for $A\in\T(X)$, the function
$\l\mapsto\a_{\l}(A)$ is norm-continuous). Moreover, from \eqref{eq:fully_coisometric_2__gauge}
and \cite[Definition 2.4]{Exel_Circle_Actions}, the $k$th spectral
subspace for $\a$ is $\T_{k}(X)$. Hence the gauge action $\a$ is
semi-saturated, because $S_{0}(\M)=\varphi_{\infty}(\M)$ and $S_{1}(X(1))$
are enough to generate $\T(X)$ as a $C^{*}$-algebra (there is a
subtlety here- see Remark \ref{rem:fully_coisometric_2}). Consequently,
by \cite[Theorem 3.1]{Abadie_Eilers_Exel}, \[
\T(X)\cong\T_{0}(X)\rtimes_{\T_{1}(X)}\Z.\]
To elaborate, if this isomorphism is implemented by $\phi$ and $\iota_{0},\iota_{1}$
are the crossed product maps from $\T_{0}(X),\T_{1}(X)$, respectively,
to $\T_{0}(X)\rtimes_{\T_{1}(X)}\Z$, then the following diagram commutes:
\begin{equation}
\xymatrix{(\T_{0}(X),\T_{1}(X))\ar[d]_{(\iota_{0},\iota_{1})}\ar[rd]\\
\T_{0}(X)\rtimes_{\T_{1}(X)}\Z\ar[r]\sb-\phi & \T(X)}
\label{eq:fully_coisometric_2__diag_T}\end{equation}
From Lemma \ref{lem:fully_coisometric_2} the pair of mappings $(\pi_{0},\pi_{1})$
is a covariant representation of $(\T_{0}(X),\T_{1}(X))$ on $\H$
in the sense of \cite[Definition 2.1]{Abadie_Eilers_Exel} (properties
(i)-(iv) therein are proved by considering first only vectors in the
total subsets of $\T_{0}(X),\T_{1}(X)$ consisting of the $S$-monomials
of degree $0,1$ respectively, and then using linearity and continuity;
notice that $\pi_{A}(\cdot)$ is missing from (iii) and (iv)). As
a result, the universality property of the crossed product implies
that there is a $C^{*}$-representation $\theta:\T_{0}(X)\rtimes_{\T_{1}(X)}\Z\to B(\H)$
such that the following diagram commutes: \[
\xymatrix{(\T_{0}(X),\T_{1}(X))\ar[d]_{(\iota_{0},\iota_{1})}\ar[rd]^{(\pi_{0},\pi_{1})}\\
\T_{0}(X)\rtimes_{\T_{1}(X)}\Z\ar[r]\sb-\theta & B(\H)}
\]
Define $\pi:=\theta\circ\phi^{-1}$. Since the diagram \eqref{eq:fully_coisometric_2__diag_T}
commutes, one infers that $\pi(A)=\pi_{k}(A)$ for $A\in\T_{k}(X)$
($k=0,1$). The proof is now complete because, as already mentioned,
$\T_{0}(X)$ and $\T_{1}(X)$ generate $\T(X)$ as a $C^{*}$-algebra.\end{proof}
\begin{rem}
\label{rem:fully_coisometric_2}If $\zeta,\eta\in X(1)$ then $S_{2}(p_{2}(\zeta\tensor\eta))=S_{1}(\zeta)S_{1}(\eta)$,
and as the span of the set of {}``simple tensors'' is \emph{norm}-dense
in $X(1)\tensor X(1)$ \emph{in the }$C^{*}$\emph{-setting}, the
$C^{*}$-algebra generated by $\T_{0}(X)$ and $\T_{1}(X)$ contains
$S_{2}(X(2))$, and similarly $S_{k}(X(k))$ for every $k\in\Z_{+}$.
Nonetheless, in the $W^{*}$-setting, this set of vectors is guaranteed
only to be $s$-dense in the tensor product (see Remark \ref{rem:self_dual_completion}).
Consequently, the $C^{*}$-algebra generated by $\T_{0}(X)$ and $\T_{1}(X)$
is ultraweakly dense in $\T(X)$, but there is no reason to expect
an equality to hold generally. We consider this limitation quite reasonable,
as the $W^{*}$-setting is more naturally suitable for ultraweakly
continuous representations of ultraweakly closed algebras.
\end{rem}
By Theorem \ref{thm:fully_coisometric_2}, condition \eqref{eq:fully_coisometric_C_star_representation}
implies $C^{*}$-extendability. While we would like to establish the
inverse implication---at least for subproduct systems with $E$ a
finite dimensional Hilbert space (analogously to Corollary \ref{cor:rel_isom_1_finite_dim_Hilbert_space})---the
fully coisometric case is more involved than the pure case. In fact,
to our knowledge, no general criterion for $C^{*}$-extendability
has yet been found even in this fundamental case. However, we will
prove that condition \eqref{eq:fully_coisometric_C_star_representation}
is indeed equivalent to $C^{*}$-extendability in the most well-known
subproduct system: the \emph{symmetric} one.
\begin{example}
\label{exa:fully_coisometric_symmetric_case}Let us consider the symmetric
subproduct system (see \cite[Example 1.3]{Subproduct_systems_2009}),
defined as follows. For a fixed $d\in\N$, take $E$ to be the Hilbert
space $\C^{d}$ (of course, $\M=\C$). Given $n\in\N$, define a projection
$p_{n}:E^{\tensor n}\to E^{\tensor n}$ by \begin{equation}
p_{n}(f_{1}\tensor\cdots\tensor f_{n}):=\frac{1}{n!}\sum_{\pi}f_{\pi(1)}\tensor\cdots\tensor f_{\pi(n)}\label{eq:symmetric_projection}\end{equation}
for all $f_{1},\ldots,f_{n}\in E$, the sum being taken over all permutations
$\pi$ of $\left\{ 1,\ldots,n\right\} $. Set $E^{\circledS n}:=p_{n}E^{\tensor n}$
(the $n$-fold \emph{symmetric tensor product} of $E$) and $E^{\circledS0}:=\C$.
The symmetric subproduct system is $\text{SSP}_{d}:=\left(E^{\circledS n}\right)_{n\in\Z_{+}}$.

If $\left\{ e_{1},\ldots,e_{d}\right\} $ is a fixed orthonormal basis
of $E$, write $T_{i}:=T_{1}(e_{i})$, $1\leq i\leq d$. There is
a bijection between all completely contractive, covariant representations
of $\text{SSP}_{d}$ on $\H$ and all commuting row contractions over
$\H$ (of length $d$) implemented by $T\mapsto(T_{1},\ldots,T_{d})$,
and we shall identify these two types of objects.

A commuting row contraction $T$ is called \emph{spherical} if $T_{i}$
is normal for all $i=1,\ldots,d$ and $T_{1}T_{1}^{*}+\ldots+T_{d}T_{d}^{*}=I_{\H}$.
A spherical row contraction is evidently fully coisometric. It was
proved in \cite[\S 8]{Arveson_subalgebras_3} that a fully coisometric,
commuting row contraction $T$ extends to a $C^{*}$-representation
if and only if $T$ is spherical.\end{example}
\begin{thm*}
Let $T$ be a fully coisometric, commuting row contraction. Then $T$
is spherical if and only if condition \eqref{eq:fully_coisometric_C_star_representation}
holds.
\end{thm*}
Sufficiency is a direct byproduct of Theorem \ref{thm:fully_coisometric_2},
but the proof will not require this fact.
\begin{proof}
Fix a fully coisometric, commuting row contraction $T$ over $\H$.
If $m\in\N$, $\eta\in E^{\circledS m}$ and $h\in\H$, we will show
that \[
\lim_{\ell\to\infty}\bigl\Vert(p_{\ell}\tensor I_{\H})(\eta\tensor\widetilde{T}_{\ell-m}^{*}h)\bigr\Vert_{E^{\circledS\ell}\tensor\H}=\norm{T_{m}(\eta)^{*}h},\]
from which the theorem's conclusion follows. Fixing an orthogonal
base $\left\{ e_{1},\ldots,e_{d}\right\} $ for $E$, it suffices
to consider only the case $\eta=e_{p}^{\tensor m}$, $1\leq p\leq d$,
because $\left\{ k^{\tensor m}:k\in E\right\} $ spans $E^{\circledS m}$,
and from \eqref{eq:fully_coisometric_C_star_representation_2}, condition
\eqref{eq:fully_coisometric_C_star_representation} being true for
all $\eta$ is equivalent to it being true for $\eta$ in a spanning
subset of $E^{\circledS m}$.

Write $[d]:=\left\{ 1,2,\ldots,d\right\} $. It follows from \cite[Proposition 6.9]{Subproduct_systems_2009}
(or a direct calculation) that for $h\in\H$ and $\ell\in\N$ we have
\begin{equation}
\widetilde{T}_{\ell}^{*}h=\sum_{(i_{1},\ldots,i_{\ell})\in[d]^{\ell}}e_{i_{1}}\tensor\cdots\tensor e_{i_{\ell}}\tensor\left(T_{i_{\ell}}^{*}\cdots T_{i_{1}}^{*}h\right).\label{eq:T_tilde_star_in_finite_d_Hilbert_spaces}\end{equation}
We introduce some notations. Given $\ell\in\N$, define an equivalence
relation over $[d]^{\ell}$ by saying that two $\ell$-tuples are
equivalent if each is a rearrangement of the other. Let $A_{\ell}$
be a subset of $[d]^{\ell}$ that contains exactly one element of
each equivalence class of $[d]^{\ell}$. If $(i_{1},\ldots,i_{\ell})\in[d]^{\ell}$
and $1\leq k\leq d$, write $c_{k}=c_{k}(i_{1},\ldots,i_{\ell})$
for the number of appearances of $k$ in the series $i_{1},\ldots,i_{\ell}$.
Denote the equivalence class of $(i_{1},\ldots,i_{\ell})$ by $[(i_{1},\ldots,i_{\ell})]_{\ell}$.
The cardinality of $[(i_{1},\ldots,i_{\ell})]_{\ell}$, denoted by
$a_{(i_{1},\ldots,i_{\ell})}$, equals $\binom{\ell}{c_{1}\,\cdots\, c_{d}}=\frac{\ell!}{c_{1}!\cdots c_{d}!}$
(the multinomial coefficient). Rewrite \eqref{eq:T_tilde_star_in_finite_d_Hilbert_spaces}
as \begin{equation}
\widetilde{T}_{\ell}^{*}h=\sum_{(j_{1},\ldots,j_{\ell})\in A_{\ell}}\Bigl(\sum_{(i_{1},\ldots,i_{\ell})\in[(j_{1},\ldots,j_{\ell})]_{\ell}}e_{i_{1}}\tensor\cdots\tensor e_{i_{\ell}}\Bigr)\tensor\left(T_{j_{\ell}}^{*}\cdots T_{j_{1}}^{*}h\right)\label{eq:T_tilde_star_symmetric}\end{equation}
(this is possible only by virtue of the commutativity of $T_{1},\ldots,T_{d}$).
As a result, \begin{equation}
\norm h^{2}=\bigl\Vert\widetilde{T}_{\ell}^{*}h\bigr\Vert^{2}=\sum_{(j_{1},\ldots,j_{\ell})\in A_{\ell}}a_{(j_{1},\ldots,j_{\ell})}\norm{T_{j_{\ell}}^{*}\cdots T_{j_{1}}^{*}h}^{2}.\label{eq:T_tilde_star_symmetric_norm_sq}\end{equation}

Given $(j_{1},\ldots,j_{\ell-m})\in[d]^{\ell-m}$, write $(\eta,j_{1},\ldots,j_{\ell-m})$
for the $\ell$-tuple $(p,\ldots,p,j_{1},\ldots,j_{\ell-m})$, and
observe that $p_{\ell}(\eta\tensor e_{i_{1}}\tensor\cdots\tensor e_{i_{\ell-m}})$
is the same for all $(i_{1},\ldots,i_{\ell-m})$ in $[(j_{1},\ldots,j_{\ell-m})]_{\ell-m}$,
and equals \[
\frac{1}{\ell!}\cdot\frac{\ell!}{a_{(\eta,j_{1},\ldots,j_{\ell-m})}}\cdot\sum_{(i_{1}.\ldots,i_{\ell})\in[(\eta,j_{1},\ldots,j_{\ell-m})]_{\ell}}e_{i_{1}}\tensor\cdots\tensor e_{i_{\ell}}\]
($\frac{\ell!}{a_{(\eta,j_{1},\ldots,j_{\ell-m})}}$ is exactly the
number of times each summand is repeated in \eqref{eq:symmetric_projection}).
Therefore, as $[(j_{1},\ldots,j_{\ell-m})]_{\ell-m}$ consists of
$a_{(j_{1},\ldots,j_{\ell-m})}$ elements, we infer from \eqref{eq:T_tilde_star_symmetric}
that \begin{multline*}
(p_{\ell}\tensor I_{\H})(\eta\tensor\widetilde{T}_{\ell-m}^{*}h)=\\
\sum_{\substack{(j_{1},\ldots,j_{\ell-m})\in\\
A_{\ell-m}}
}\frac{a_{(j_{1},\ldots,j_{\ell-m})}}{a_{(\eta,j_{1},\ldots,j_{\ell-m})}}\Bigl(\sum_{\substack{(i_{1}.\ldots,i_{\ell})\in\\
{}[(\eta,j_{1},\ldots,j_{\ell-m})]_{\ell}}
}e_{i_{1}}\tensor\cdots\tensor e_{i_{\ell}}\Bigr)\tensor\left(T_{j_{\ell-m}}^{*}\cdots T_{j_{1}}^{*}h\right).\end{multline*}
Write $c_{p}=c_{p}(j_{1},\ldots,j_{\ell-m})$. All summands in the
last sum are mutually orthogonal, hence we have \begin{multline}
\bigl\Vert(p_{\ell}\tensor I_{\H})(\eta\tensor\widetilde{T}_{\ell-m}^{*}h)\bigr\Vert^{2}=\\
\begin{split} & =\sum_{(j_{1},\ldots,j_{\ell-m})\in A_{\ell-m}}\frac{a_{(j_{1},\ldots,j_{\ell-m})}^{2}}{a_{(\eta,j_{1},\ldots,j_{\ell-m})}^{2}}a_{(\eta,j_{1},\ldots,j_{\ell-m})}\bigl\Vert T_{j_{\ell-m}}^{*}\cdots T_{j_{1}}^{*}h\bigr\Vert^{2}\\
 & =\sum_{\text{the same}}\frac{(c_{p}+m)\cdots(c_{p}+1)}{\ell(\ell-1)\cdots(\ell-m+1)}a_{(j_{1},\ldots,j_{\ell-m})}\bigl\Vert T_{j_{\ell-m}}^{*}\cdots T_{j_{1}}^{*}h\bigr\Vert^{2}\end{split}
\label{eq:symmetric_example_1}\end{multline}
(because $c_{k}(\eta,j_{1},\ldots,j_{\ell-m})=c_{k}$ for all $k\neq p$
and $c_{p}(\eta,j_{1},\ldots,j_{\ell-m})=c_{p}+m$). Next, we assert
that the sum of all summands in \eqref{eq:symmetric_example_1} corresponding
to tuples $(j_{1},\ldots,j_{\ell-m})\in A_{\ell-m}$ with $c_{p}<\sqrt{\ell}$
are negligible as $\ell\to\infty$. For, summing over these tuples
alone, the result is dominated by \begin{multline}
\left(\frac{\sqrt{\ell}+m}{\ell-m+1}\right)^{m}\sum_{\substack{\text{such tuples}\\
(j_{1},\ldots,j_{\ell-m})}
}a_{(j_{1},\ldots,j_{\ell-m})}\bigl\Vert T_{j_{\ell-m}}^{*}\cdots T_{j_{1}}^{*}h\bigr\Vert^{2}\leq\\
\left(\frac{\sqrt{\ell}+m}{\ell-m+1}\right)^{m}\sum_{(j_{1},\ldots,j_{\ell-m})\in A_{\ell-m}}a_{(j_{1},\ldots,j_{\ell-m})}\bigl\Vert T_{j_{\ell-m}}^{*}\cdots T_{j_{1}}^{*}h\bigr\Vert^{2}=\left(\frac{\sqrt{\ell}+m}{\ell-m+1}\right)^{m}\norm h^{2}\label{eq:symmetric_example_diff_1}\end{multline}
(see \eqref{eq:T_tilde_star_symmetric_norm_sq}) and the right side
converges to $0$ as $\ell\to\infty$. Consequently, letting $\ell$
grow big enough, we may assume (in particular) that $c_{p}\geq m$.
For such $(\ell-m)$-tuple $(j_{1},\ldots,j_{\ell-m})$, write $(k_{1},\ldots,k_{\ell-2m})$
for an $(\ell-2m)$-tuple obtained by removing (any) $m$ repetitions
of $p$ from $(j_{1},\ldots,j_{\ell-m})$. So in conclusion, after
removing all $(j_{1},\ldots,j_{\ell-m})\in A_{\ell-m}$ with $c_{p}<\sqrt{\ell}$
from the sum in \eqref{eq:symmetric_example_1}, we are left with
\begin{multline}
\sum_{\substack{(j_{1},\ldots,j_{\ell-m})\in A_{\ell-m}\\
\text{with }c_{p}\ge\sqrt{\ell}}
}\frac{(c_{p}+m)\cdots(c_{p}+1)}{\ell(\ell-1)\cdots(\ell-m+1)}a_{(j_{1},\ldots,j_{\ell-m})}\bigl\Vert T_{j_{\ell-m}}^{*}\cdots T_{j_{1}}^{*}h\bigr\Vert^{2}=\\
\begin{split}=\sum_{\text{the same}} & \frac{(c_{p}+m)\cdots(c_{p}+1)}{\ell(\ell-1)\cdots(\ell-m+1)}\frac{(\ell-m)\cdots(\ell-2m+1)}{c_{p}\cdots(c_{p}-m+1)}\cdot\\
 & \cdot a_{(k_{1},\ldots,k_{\ell-2m})}\bigl\Vert T_{k_{\ell-2m}}^{*}\cdots T_{k_{1}}^{*}T_{p}^{m*}h\bigr\Vert^{2}.\end{split}
\label{eq:symmetric_example_2}\end{multline}
Since $m$ is fixed, $\frac{(\ell-m)\cdots(\ell-2m+1)}{\ell(\ell-1)\cdots(\ell-m+1)}\to1$
as $\ell\to\infty$. So for the purpose of the limit of $\bigl\Vert(p_{\ell}\tensor I_{\H})(\eta\tensor\widetilde{T}_{\ell-m}^{*}h)\bigr\Vert^{2}$
we can remove this factor from \eqref{eq:symmetric_example_2}. Additionally,
since $c_{p}\ge\sqrt{\ell}$, $\frac{(c_{p}+m)\cdots(c_{p}+1)}{c_{p}\cdots(c_{p}-m+1)}$
converges to $1$ as $\ell\to\infty$ \emph{uniformly} for all relevant
tuples. Removing this factor as well from \eqref{eq:symmetric_example_2}
we get \begin{multline*}
\sum_{\substack{(j_{1},\ldots,j_{\ell-m})\in A_{\ell-m}\\
\text{with }c_{p}\ge\sqrt{\ell}}
}a_{(k_{1},\ldots,k_{\ell-2m})}\bigl\Vert T_{k_{\ell-2m}}^{*}\cdots T_{k_{1}}^{*}T_{p}^{m*}h\bigr\Vert^{2}\\
=\sum_{\substack{(i_{1},\ldots,i_{\ell-2m})\in A_{\ell-2m}\\
\text{with }c_{p}(i_{1},\ldots,i_{\ell-2m})\ge\sqrt{\ell}-m}
}a_{(i_{1},\ldots,i_{\ell-2m})}\bigl\Vert T_{i_{\ell-2m}}^{*}\cdots T_{i_{1}}^{*}T_{p}^{m*}h\bigr\Vert^{2}.\end{multline*}
Had the $c_{p}$ of the $(\ell-2m)$-tuples in the last sum not been
restricted, we would have had the sum \[
\sum_{(i_{1},\ldots,i_{\ell-2m})\in A_{\ell-2m}}a_{(i_{1},\ldots,i_{\ell-2m})}\bigl\Vert T_{i_{\ell-2m}}^{*}\cdots T_{i_{1}}^{*}T_{p}^{m*}h\bigr\Vert^{2}=\bigl\Vert T_{p}^{m*}h\bigr\Vert^{2}=\bigl\Vert T_{m}(\eta)^{*}h\bigr\Vert^{2},\]
(by \eqref{eq:T_tilde_star_symmetric_norm_sq}), as desired. So all
that is left is to see that the difference between the two sums converges
to zero. This will be done exactly as before (see \eqref{eq:symmetric_example_diff_1}),
denoting now $c_{p}=c_{p}(i_{1},\ldots,i_{\ell-2m})$: \begin{multline*}
\sum_{\substack{(i_{1},\ldots,i_{\ell-2m})\in A_{\ell-2m}\\
\text{with }c_{p}<\sqrt{\ell}-m}
}a_{(i_{1},\ldots,i_{\ell-2m})}\bigl\Vert T_{i_{\ell-2m}}^{*}\cdots T_{i_{1}}^{*}T_{p}^{m*}h\bigr\Vert^{2}\\
\begin{split} & =\sum_{\text{the same}}\frac{(c_{p}+m)\cdots(c_{p}+1)}{(\ell-m)\cdots(\ell-2m+1)}a_{(\eta,i_{1},\ldots,i_{\ell-2m})}\bigl\Vert T_{i_{\ell-2m}}^{*}\cdots T_{i_{1}}^{*}T_{p}^{m*}h\bigr\Vert^{2}\\
 & \leq\left(\frac{\sqrt{\ell}}{\ell-2m+1}\right)^{m}\sum_{\text{the same}}a_{(\eta,i_{1},\ldots,i_{\ell-2m})}\bigl\Vert T_{i_{\ell-2m}}^{*}\cdots T_{i_{1}}^{*}T_{p}^{m*}h\bigr\Vert^{2}\\
 & \leq\left(\frac{\sqrt{\ell}}{\ell-2m+1}\right)^{m}\sum_{(t_{1},\ldots,t_{\ell-m})\in A_{\ell-m}}a_{(t_{1},\ldots,t_{\ell-m})}\bigl\Vert T_{t_{\ell-m}}^{*}\cdots T_{t_{1}}^{*}h\bigr\Vert^{2}\\
 & =\left(\frac{\sqrt{\ell}}{\ell-2m+1}\right)^{m}\norm h^{2}\xrightarrow[\ell\to\infty]{}0\end{split}
\end{multline*}
(see \eqref{eq:T_tilde_star_symmetric_norm_sq}).
\end{proof}

\section{Conclusions and examples\label{sec:Conclusions-and-examples}}

In this section we combine some of the previous results, and give
more examples. We begin with a general Wold decomposition. Recall
that $Q=\slim_{n\to\infty}\widetilde{T}_{n}\widetilde{T}_{n}^{*}$.
\begin{thm}
\label{thm:rel_isom_and_fully_coiso}Let $X=\left(X(n)\right)_{n\in\Z_{+}}$
be a standard subproduct system in the $C^{*}$-setting, and suppose
that $T$ is a completely contractive, covariant representation of
$X$ on $\H$. Assume that $T$ is relatively isometric and \begin{equation}
\lim_{\ell\to\infty}\bigl\Vert(p_{\ell}\tensor Q)(\eta\tensor\widetilde{T}_{\ell-m}^{*}h)\bigr\Vert_{X(\ell)\tensor\H}=\norm{T_{m}(\eta)Qh}_{\H}\label{eq:rel_isom_and_fully_coiso}\end{equation}
for all $m\in\N$, $\eta\in X(m)$ and $h\in\H$. Then there exist
a Hilbert space $\U$, a unitary operator $W:\H\to(\F_{X}\tensor_{\sigma}\D)\oplus\U$
and a fully coisometric, covariant representation $Z=\left(Z_{n}\right)_{n\in\Z_{+}}$
of $X$ on $\U$, which extends to a $C^{*}$-representation, such
that $WT_{n}(\zeta)W^{-1}=\left(S_{n}(\zeta)\tensor I_{\D}\right)\oplus Z_{n}(\zeta)$
for all $n\in\Z_{+}$, $\zeta\in X(n)$. In particular, $T$ extends
to a $C^{*}$-representation.\end{thm}
\begin{proof}
From Theorems \ref{thm:dilation_1} and \ref{thm:rel_isom_1} and
their proofs, we should ascertain that the covariant representation
$Z$ (of $X$ on the Hilbert subspace $\U$ of $\H$) mentioned there
extends to a $C^{*}$-representation. Under present assumptions, $Q$
is the projection of $\H$ on $\U$ and $Y$ is $Q$ with codomain
$\U$. In particular, $\U=\Img Q=\Img Y$. By \eqref{eq:dilation_thm_Y_adj_Z_tilde},
\[
(p_{\ell}\tensor I_{\U})(\eta\tensor\widetilde{Z}_{\ell-m}^{*}Yh)=(p_{\ell}\tensor I_{\U})\bigl(\eta\tensor(I_{X(\ell-m)}\tensor Y)\widetilde{T}_{\ell-m}^{*}h\bigr)=(p_{\ell}\tensor Q)(\eta\tensor\widetilde{T}_{\ell-m}^{*}h).\]
Since $Z_{m}(\eta)Yh=T_{m}(\eta)Qh$, we infer from \eqref{thm:rel_isom_and_fully_coiso}
that $Z$ satisfies the conditions of Theorem \ref{thm:fully_coisometric_2}.
This completes the proof.\end{proof}
\begin{example}
[Homogeneous ideals] We extend results of \cite[\S 2]{Quotient_Tensor_Algebras},
showing that most of their hypothesis (H1) is superfluous (also cf.~\cite{Sarason_interpolation},
\cite[\S 2]{Davidson_Pitts_Neva_1998}, \cite[\S 1]{Arias_Popescu_I}
and \cite[\S 7]{Subproduct_systems_2009}). Throughout this example
we shall be working in the $C^{*}$-setting. Let $E$ be a $C^{*}$-correspondence
that is essential as a left $\M$-module. Write $S=(S_{n})_{n\in\Z_{+}}$
for the shift in $\F(E)$. Given a two-sided, norm closed ideal $\I\trianglelefteq\T_{+}(E)$,
define $\I^{(n)}:=\I\cap\left\{ S_{n}(\zeta):\zeta\in E^{\tensor n}\right\} $
for $n\in\Z_{+}$. We call $\I$ \emph{homogeneous} if $\I=\overline{\linspan}\bigcup_{n\in\Z_{+}}\I^{(n)}$
(equivalently: $\I$ is invariant under the gauge action on $\T_{+}(E)$).
Our standing hypotheses are:
\begin{enumerate}
\item \label{enu:homog_ideal_1}The ideal $\I\trianglelefteq\T_{+}(E)$
is homogeneous.
\item \label{enu:homog_ideal_2}The submodule $Y(n):=\left\{ \zeta\in E^{\tensor n}:S_{n}(\zeta)\in\I\right\} $
is orthogonally complementable in $E^{\tensor n}$ for all $n$ and
$Y(0)=\left\{ 0\right\} $.
\end{enumerate}
Given $n\in\Z_{+}$, define $X(n)$ to be the orthogonal complement
of $Y(n)$ in $E^{\tensor n}$, and let $p_{n}$ denote the orthogonal
projection of $E^{\tensor n}$ on $X(n)$. If $n,m\in\Z_{+}$, then
since $\I$ is an ideal we obtain $Y(n)\tensor E^{\tensor m}\subseteq Y(n+m)$
and $E^{\tensor n}\tensor Y(m)\subseteq Y(n+m)$, that is, $p_{n+m}\leq p_{n}\tensor I_{E^{\tensor m}}$
and $p_{n+m}\leq I_{E^{\tensor n}}\tensor p_{m}$. Therefore, since
$p_{n}\tensor I_{E^{\tensor m}}$ and $I_{E^{\tensor n}}\tensor p_{m}$
commute, $X(n)\tensor X(m)=\Img(p_{n}\tensor p_{m})=\bigl(\Img(p_{n}\tensor I_{E^{\tensor m}})\bigr)\cap\bigl(\Img(I_{E^{\tensor n}}\tensor p_{m})\bigr)\supseteq\Img p_{n+m}=X(n+m)$.
In conclusion, $X=\left(X(n)\right)_{n\in\Z_{+}}$ is a subproduct
system with $X(0)=\M$. Moreover, \ref{enu:homog_ideal_1} yields
that \begin{multline*}
\overline{\I\F(E)}=\overline{\bigl(\overline{\linspan}\bigcup_{n\in\Z_{+}}\I^{(n)}\bigr)\F(E)}=\overline{\bigl(\linspan\bigcup_{n\in\Z_{+}}\I^{(n)}\bigr)\F(E)}\\
=\overline{\linspan\bigcup_{n,m\in\Z_{+}}\I^{(n)}E^{\tensor m}}=\bigoplus_{k\in\Z_{+}}\overline{\linspan\bigcup_{\substack{n,m\in\Z_{+}\\
n+m=k}
}\I^{(n)}E^{\tensor m}}=\bigoplus_{k\in\Z_{+}}Y(k)\end{multline*}
(for $\I^{(n)}E^{\tensor m}\subseteq Y(n)\tensor E^{\tensor m}\subseteq Y(n+m)=\overline{\I^{(n+m)}\M}$).
Consequently, \begin{equation}
\F_{X}=\F(E)\ominus\overline{\I\F(E)}.\label{eq:homog_ideal_F_X}\end{equation}
Let $P:=\bigoplus_{n\in\Z_{+}}p_{n}$ be the projection of $\F(E)$
on $\F_{X}$. The relation between $\T_{+}(E)$ and $\T_{+}(X)$ is
illuminated in the following result.
\begin{thm*}
$\T_{+}(X)\cong\T_{+}(E)\slash\I$ in the sense that the operator
algebras are (canonically) completely isometrically isomorphic.\end{thm*}
\begin{proof}
Write $S^{X}=(S_{n}^{X})_{n\in\Z_{+}}$ for the $X$-shift. Define
$\kappa:\T_{+}(E)\slash\I\to\T_{+}(X)$ by $T+\I\mapsto PT_{|\F_{X}}$
for $T\in\T_{+}(E)$. The mapping $\kappa$ is well-defined by \eqref{eq:homog_ideal_F_X}
and $\kappa(S_{n}(\zeta)+\I)=S_{n}^{X}(\zeta)$ for all $n\in\Z_{+}$,
$\zeta\in X(n)$. It is contractive since if $T\in\T_{+}(E)$ and
$Z\in\I$, then from \eqref{eq:homog_ideal_F_X}, \[
\bigl\Vert PT_{|\F_{X}}\bigr\Vert=\bigl\Vert PTP\bigr\Vert=\bigl\Vert P(T+Z)P\bigr\Vert\leq\bigl\Vert T+Z\bigr\Vert,\]
and so $\bigl\Vert PT_{|\F_{X}}\bigr\Vert\leq\inf_{Z\in\I}\bigl\Vert T+Z\bigr\Vert=\bigl\Vert T+\I\bigr\Vert$.
Using the identification of $M_{n}(\T_{+}(E)\slash\I)$ with $M_{n}(\T_{+}(E))\slash M_{n}(\I)$,
one proves in a similar fashion that $\kappa$ is actually completely
contractive.

Let $\pi$ be a nondegenerate completely isometric representation
of $\T_{+}(E)\slash\I$ on some Hilbert space $\H$ so that $\M\ni a\mapsto\pi(\varphi_{\infty}(a)+\I)$
is a (nondegenerate, because $E$ is essential) $C^{*}$-representation
of $\M$. Define a completely contractive, covariant representation
$T$ of $X$ on $\H$ by $T_{n}(\zeta):=\pi(S_{n}(\zeta)+\I)$ ($n\in\Z_{+},\zeta\in X(n)$).
Indeed, if $n,m\in\Z_{+}$, $\zeta\in X(n)$ and $\eta\in X(m)$,
we obtain \[
T_{n}(\zeta)T_{m}(\eta)=\pi(S_{n}(\zeta)S_{m}(\eta)+\I)=\pi(S_{n+m}(\zeta\tensor\eta)+\I)=\pi(S_{n+m}(p_{n+m}(\zeta\tensor\eta))+\I)\]
(as $(I_{E^{\tensor(n+m)}}-p_{n+m})(\zeta\tensor\eta)\in E^{\tensor(n+m)}\ominus X(n+m)$).
Moreover, since $S_{n}(\cdot)$ is completely contractive for all
$n$, the same is true for $S_{n}(\cdot)+\I$, and hence also for
$T_{n}(\cdot)$. Therefore, Theorem \ref{thm:c_c_c_representation__to__A_morphism}
implies that there is a completely contractive linear mapping $\vartheta:\T_{+}(X)\to\T_{+}(E)\slash\I$
satisfying $\vartheta(S_{n}^{X}(\zeta))=S_{n}(\zeta)+\I$ for all
$n\in\Z_{+}$, $\zeta\in X(n)$. Hence $\vartheta=\kappa^{-1}$. Both
$\kappa$ and $\vartheta$ are completely contractive, thus they are
completely isometric. This completes the proof.\end{proof}
\begin{cor*}
$d(T,\I)=\norm{PTP}$ for every $T\in\T_{+}(E)$.
\end{cor*}
\end{example}

\begin{example}
Fix a von Neumann algebra $\M\subseteq B(\H)$ and a $cp$-semigroup
$\left(\Theta_{n}\right)_{n\in\Z_{+}}$ on $\M$. It is established
in \cite[Theorem 2.2]{Subproduct_systems_2009} and the preceding
discussion that there exist a subproduct system $X=\left(X(n)\right)_{n\in\Z_{+}}$
over the commutant $\M'$ and a completely contractive, covariant
representation $T=\left(T_{n}\right)_{n\in\Z_{+}}$ of $X$ on $\H$
such that \begin{equation}
(\A n\in\Z_{+},a\in\M)\quad\quad\Theta_{n}(a)=\widetilde{T}_{n}(I_{X(n)}\tensor a)\widetilde{T}_{n}^{*}.\label{eq:cp_semigroup_example}\end{equation}
This construction is interesting in view of \cite[Theorem 2.6]{Subproduct_systems_2009}.
Let us see when this $T$ is relatively isometric. From \eqref{eq:cp_semigroup_example},
$T_{n}$ is a partial isometry if and only if \emph{$\Theta_{n}(I)$
}is a projection. Using the notations of \cite{Subproduct_systems_2009},
we have $\Delta_{*}=I-\Theta_{1}(I)$ and $T_{n}(x)=W_{\Theta_{n}}^{*}x$
for $x\in X(n)=\L_{\M}(\H,\M\tensor_{\Theta_{n}}\H)$. Consequently,
\eqref{eq:relatively_isometric_operators} is equivalent to \[
(\A n\in\Z_{+},x\in X(n))\quad\quad(I-\Theta_{1}(I))x^{*}W_{\Theta_{n}}W_{\Theta_{n}}^{*}x(I-\Theta_{1}(I))=x^{*}x(I-\Theta_{1}(I))\]
($\sigma$ is the identity representation of $\M'$). For instance,
this condition is easily satisfied when the operators $W_{\Theta_{n}}$
are all coisometric, which happens, e.g., when $\Theta$ is an $e$-semigroup
(this in turn implies that $T$ is actually \emph{isometric}, so that
$X$ is necessarily a \emph{product} system by \cite[Corollary 2.8]{Subproduct_systems_2009}).
\end{example}

\begin{example}
[Strictly cyclic subproduct systems]Let $\M$ be a $C^{*}$-algebra.
A $C^{*}$-correspondence $F$ over $\M$ is\emph{ }called strictly
cyclic\emph{ }(\cite[p.~419]{Tensor_algebras}) if $F=P\M$ for some
projection $P$ in $M(\M)$, the multiplier algebra of $\M$ (the
left action of $\M$ on $F$ is not restricted). A subproduct system
$X=(X(n))_{n\in\Z_{+}}$ over $\M$ is strictly cyclic if $E=X(1)=\M$
and $X(n)$ is strictly cyclic---i.e., $X(n)=P_{n}\M^{\tensor n}$---for
all $n\geq2$. When given such $X$, the projection $p_{n}$ of $\M^{\tensor n}$
on $X(n)$ is of a very concrete form: it is multiplication on the
left by $P_{n}$, resulting in a more simplified condition \eqref{eq:fully_coisometric_C_star_representation}.

For example, define the left operation of $\M$ on $E:=\M$ to be
left multiplication. Then $E^{\tensor n}$ may be identified with
$\M$ via $a_{1}\tensor\cdots\tensor a_{n}\mapsto a_{1}\cdots a_{n}$,
and if the projections $P_{n}\in M(\M)$, $n\geq2$, commute with
$\M$ and satisfy $P_{n+m}\M\subseteq P_{n}P_{m}\M$ for each $n,m\in\N$
($P_{1}:=I$), we obtain a strictly cyclic subproduct system by setting
$X(n):=P_{n}\M$.
\end{example}

\section*{Acknowledgments}

\noindent The author would like to express his deep gratitude to Baruch
Solel for many intriguing and fruitful conversations on the content
of this paper. He is grateful to the referee for carefully reading
the manuscript and making useful suggestions.

\end{document}